\documentclass[11pt,a4paper,reqno]{amsart}

\usepackage[a4paper,
            textwidth=15.5cm,
            textheight=23cm,
            hcentering,
            vcentering,
            headheight=14pt,
            headsep=0.6cm,
            footskip=1.1cm]{geometry}

\usepackage[T1]{fontenc}
\usepackage[utf8]{inputenc}
\usepackage{amsmath,amsthm,amscd}
\usepackage{mathtools}
\IfFileExists{newpxmath.sty}{\usepackage[scaled=0.97]{newpxtext}   
\usepackage[scaled=0.97]{newpxmath}}{\usepackage{amssymb,mathpazo}}   
\usepackage{microtype}
\IfFileExists{stmaryrd.sty}{\usepackage{stmaryrd}}{}                 

\usepackage{array,booktabs,mathrsfs,graphicx}
\usepackage{caption}
\usepackage{verbatim}
\usepackage{enumitem}
\usepackage{xcolor}
\definecolor{linkblue}{RGB}{0,51,128}
\usepackage[colorlinks=true,
            linkcolor=linkblue,
            citecolor=linkblue,
            urlcolor=linkblue,
            pdfusetitle,
            bookmarksnumbered]{hyperref}
\usepackage[nameinlink,capitalize,noabbrev]{cleveref}

\allowdisplaybreaks[1]
\numberwithin{equation}{section}

\theoremstyle{plain}
\newtheorem{theorem}{Theorem}[section]
\newtheorem{proposition}[theorem]{Proposition}
\newtheorem{lemma}[theorem]{Lemma}
\newtheorem{corollary}[theorem]{Corollary}
\newtheorem*{liftcriterion}{Lemma (lifting a column)}
\newtheorem*{correctionlemma}{Lemma (endpoint corrections)}
\newtheorem*{basislemma}{Lemma (descent of a basis)}
\newtheorem{mainresult}{Theorem}
\newtheorem{maincorollary}[mainresult]{Corollary}

\newcommand{\OO}{\mathcal O}
\newcommand{\PP}{\mathbb P}
\newcommand{\Dcal}{\mathcal D}
\newcommand{\Bcal}{\mathcal B}
\newcommand{\Ccal}{\mathcal C}
\newcommand{\m}{\mathfrak m}
\DeclareMathOperator{\res}{res}
\DeclareMathOperator{\divi}{div}
\DeclareMathOperator{\Pic}{Pic}
\DeclareMathOperator{\Hilb}{Hilb}
\DeclareMathOperator{\Spec}{Spec}
\DeclareMathOperator{\im}{im}
\newcommand{\eHK}{\operatorname{e}_{\mathrm{HK}}}
\DeclareMathOperator{\rk}{rank}
\hypersetup{pdftitle={Monsky's nodal cubic and the arithmetic of Frobenius invariants},pdfauthor={Hoang Le Truong},bookmarksnumbered=true}
\setlist[enumerate]{label=(\arabic*),leftmargin=2em,itemsep=2pt,topsep=5pt}
\title[Monsky's nodal cubic and Frobenius invariants]{Monsky's nodal cubic and the arithmetic of Frobenius invariants}
\author[H.L.Truong]{Hoang Le Truong}
\address{Institute for Artificial Intelligence, VNU University of Engineering and Technology, Hanoi, Vietnam}
\email{hltruong@vnu.edu.vn, truonghoangle@gmail.com}
\subjclass[2020]{Primary 13A35. Secondary 13D40, 13H15, 14H60}
\keywords{Hilbert--Kunz multiplicity, $F$-signature, nodal cubic, syzygy bundle, Frobenius splitting series, transcendence}
\hypersetup{pdfauthor={Hoang Le Truong}}
\thanks{ChatGPT (OpenAI) was used in an iterative drafting process based on the authors' research notes and detailed instructions, and to assist with subsequent revisions of the manuscript. The authors reviewed and substantially revised the resulting text, verified the mathematical arguments and references, and take responsibility for the final manuscript.}
\date{}
\begin{document}

\begin{abstract}
We prove Monsky's conjecture on the colengths of the ideals
\((x^q,y^q,z^q,h^j)\) in \(\mathbb{F}_2[x,y,z]\), where
\(h=x^3+y^3+xyz\), \(q=2^e\), and \(0\le j\le q\).
We derive Monsky's recursion by interchanging the values of sections
of syzygy bundles on the two branches at the node of the cubic \(h=0\).
For the hypersurface \(uv+h=0\), the Hilbert--Kunz and Frobenius
splitting series are algebraic of degree two over \(\mathbb{Q}(w)\).
For \(uv+h(X_1)+\cdots+h(X_m)=0\), where the \(X_i\) are disjoint
triples of variables and \(m\ge 2\), both series are transcendental
over \(\mathbb{Q}(w)\). When \(m=2\), we express the Hilbert--Kunz
multiplicity in terms of a period of an elliptic curve in
characteristic zero. This gives an explicit seven-dimensional
hypersurface over \(\mathbb{F}_2\) with transcendental Hilbert--Kunz
multiplicity and transcendental \(F\)-signature.
\end{abstract}
\maketitle

\section{Introduction}

Let $(A,\m)$ be a Noetherian local ring of dimension $d$ and characteristic $p>0$. For $q=p^e$, the Frobenius power $\m^{[q]}$ is the ideal generated by the $q$th powers of elements of $\m$. In a regular local ring the length $\ell_A(A/\m^{[q]})$ is exactly $q^d$ \cite{Kunz}. In general the quotient of this length by $q^d$ has a limit, the Hilbert--Kunz multiplicity $\eHK(A)$ \cite{Mon83}. This limit can be irrational even when $A$ is a domain \cite{Bren13}. In this paper we prove a conjecture of Monsky on a nodal plane cubic and use it to give an explicit hypersurface over $\mathbb F_2$ with transcendental Hilbert--Kunz multiplicity. The same equation exhibits a change in arithmetic. Suspending one copy of the cubic gives a quadratic irrational, whereas suspending the sum of two copies gives a transcendental number.

Let $k$ be a field of characteristic two, let $h=x^3+y^3+xyz$, and let $C\subset\PP^2_k$ be the cubic defined by $h=0$. Its $j$th thickening is defined by $h^j=0$. For $q=2^e$ and $j\ge0$, write
\[
 E_q(j)=\dim_k k[x,y,z]/(x^q,y^q,z^q,h^j).
\]
The formula for $E_q(1)$, the colength on the reduced cubic, goes back to Pardue and was proved by Buchweitz and Chen \cite[Theorem~3]{BC}, \cite{M11}. The higher thickenings enter when one adds two variables and passes to the suspension
\[
 R=\mathbb F_2[[u,v,x,y,z]]/(uv+h).
\]
Indeed, writing $\m_R$ for its maximal ideal and taking $k=\mathbb F_2$, Monsky obtained
\[
 \ell_R(R/\m_R^{[q]})=q^3+2\sum_{j=1}^{q-1}E_q(j).
\]
He conjectured a formula for the individual summands in terms of the dyadic parameter $t=j/q$ \cite[Theorem~2.2 and Conjecture~1.5]{Mon08}, \cite[Conjecture~2.3]{M09}. In its recursive form the conjecture uses functions $\phi_m$, indexed by $m\ge0$, on the dyadic rationals in $[0,1]$. They vanish at both endpoints. The two half-interval rules in \eqref{rules} determine their values successively in order of increasing denominator. Only the initial function $\phi_0$ appears in the colength formula.

\begin{mainresult}[Monsky's nodal cubic conjecture]\label{main}
Let $k$ be a field of characteristic two. For every $q=2^e$ and every integer $0\le j\le q$,
\[
 E_q(j)=q^3\bigl(t+\phi_0(t)\bigr),\qquad t=j/q,
\]
where the functions $\phi_m$ are defined by \eqref{rules}. For $j\ge q$ one has $E_q(j)=q^3$.
\end{mainresult}

The proof recovers the two half-interval rules from the geometry at the node of $C$. Frobenius pullbacks of syzygy bundles have been used to compute Hilbert--Kunz multiplicities in graded dimension two \cite{BrenRat,TriSemi}. On the reduced nodal cubic, such a bundle can be described by its pullback to the normalization $\PP^1$ and an identification of the two fibers above the node \cite{DG,BBDG}. Burban described Frobenius pullback in these terms \cite{Burban}, and Monsky used this description to compute on the reduced cubic \cite{M11}.

On the $j$th thickening, a section on one branch must agree with a section on the other along their entire intersection. In formal parameters $a,b$ for which the cubic has local equation $ab=0$, that intersection has ring $k[[a,b]]/(a^j,b^j)$ and length $j^2$. This length accounts for the quadratic correction in the recursion. To compare the spaces of sections, we change the twist at a flex of $C$ so as to interchange the two restrictions at the node (Lemma~\ref{interchange}). The interchange allows a section to extend through the summands on the opposite branch in reverse order. Section~\ref{sec:cubic} derives the resulting colength comparison. Section~\ref{sec:colon} tracks the lengths of the chains of syzygy summands under the colon operations arising in that comparison. These changes account for the shifts of the index $m$ in \eqref{rules}.

Section~\ref{sec:proof} combines the comparison and the index shifts to prove Theorem~\ref{main}. The same formula also answers two questions of Blickle, Schwede and Tucker for this cubic \cite[Questions~4.7--4.8]{BST}. Indeed, for $B=\mathbb F_2[[x,y,z]]$ and dyadic $t=j/q$ in $[0,1]$, the $F$-signature of the pair $(B,h^t)$ is $1-q^{-3}E_q(j)$. The recursion yields rational values at every rational exponent, but no polynomial agrees with the signature function on a left neighborhood of its threshold.

We now return to the suspension. To record the effect of Frobenius at each iterate, we use two generating series. For an $F$-finite local ring $(A,\m)$ with perfect residue field, let $a_e(A)$ be the free rank of the module $F_*^e A$ obtained by restriction of scalars along the $e$th Frobenius map. The limit of $a_e(A)/p^{e\dim A}$ is the $F$-signature $s(A)$ \cite{HL02,SVdB,Tuc12}. With $w$ an indeterminate, the Hilbert--Kunz series and the Frobenius splitting series are respectively
\[
 H_A(w)=\sum_{e\ge0}\ell_A(A/\m^{[p^e]})w^e,
 \qquad S_A(w)=\sum_{e\ge0}a_e(A)w^e.
\]
For $A=R$, Monsky computed the Hilbert--Kunz series and its limiting multiplicity assuming his conjecture \cite{Mon08,Mtrans}, and Tucker deduced the corresponding signature \cite[Proposition~4.22]{Tuc12}. Theorem~\ref{main} makes these computations unconditional (Section~\ref{sec:area}).

\begin{maincorollary}\label{main:irr}
The cubic suspension $R$ has
\[
 \eHK(R)=\frac43+\frac{5}{14\sqrt7},\qquad
 s(R)=\frac23-\frac{5}{14\sqrt7}.
\]
\end{maincorollary}

\begin{maincorollary}\label{main:series}
For $R$, one has $S_R(w)=2/(1-16w)-H_R(w)$. Each series generates the quadratic extension $\mathbb Q(w,\sqrt{1-4w^2})$ of $\mathbb Q(w)$.
\end{maincorollary}

Monsky also showed how the colengths of one cubic could detect transcendence among the multiplicities of six related hypersurfaces \cite[Theorems~1.10 and~2.2]{Mtrans}. Let $X$ and $X'$ denote disjoint triples of variables, and write $h(X)$ and $h(X')$ for the same cubic evaluated in those triples. With a further variable $t$, Monsky's six rings are
\[
 A_{i,b}=\mathbb F_2[[u,v,X,X',t]]/
 (uv+h(X)+h(X')^{2^i}+t^b),
 \qquad i\in\{0,1,2\},\quad b\in\{3,5\}.
\]
His argument expresses a normalized transcendental elliptic period as a rational linear combination of their multiplicities. Applying Theorem~\ref{main} and the identity $s(A_{i,b})=2-\eHK(A_{i,b})$ gives the following result in Section~\ref{sec:period}.

\begin{mainresult}\label{main:trans}
All six rings $A_{i,b}$ are eight-dimensional strongly $F$-regular Gorenstein hypersurfaces. For at least one pair $(i,b)$, both $\eHK(A_{i,b})$ and $s(A_{i,b})$ are transcendental over $\mathbb Q$.
\end{mainresult}

Das, Landsittel and Pham established transcendence for a homogeneous primary ideal in a normal standard graded domain over every uncountable algebraically closed field of characteristic $p>2$ \cite{DLP}. Here the ideal is the maximal ideal and the ground field is $\mathbb F_2$. To identify a particular hypersurface with transcendental multiplicity, we study sums of copies of $h$ directly. For $m\ge1$, take disjoint triples $X_1,\ldots,X_m$ and put
\[
 R_m=\mathbb F_2[[u,v,X_1,\ldots,X_m]]/
 \left(uv+\sum_{i=1}^m h(X_i)\right).
\]
Thus $R_1$ is the suspension $R$. Passing from one copy to two already changes the arithmetic of both generating series.

\begin{mainresult}\label{main:repeated}
For every $m\ge2$, both $H_{R_m}$ and $S_{R_m}$ are transcendental over $\mathbb Q(w)$. Each satisfies a nonzero linear differential equation over $\mathbb Q[w]$. The ring $R_m$ is strongly $F$-regular and $s(R_m)\ge1/3$.
\end{mainresult}

Section~\ref{sec:repeated} proves the series statement by taking coefficientwise products. To pass from transcendence of a series to transcendence of the multiplicity of $R_2$, one must also control a particular value. Section~\ref{sec:r2} expresses the relevant value at $w=1/128$ in terms of an elliptic period in characteristic zero. The residues of the differential are an algebraic multiple of an integral divisor whose class has infinite order. This is the input for W\"ustholz's period theorem and yields the following conclusion.

\begin{mainresult}\label{main:r2}
The seven-dimensional hypersurface $R_2$ has transcendental Hilbert--Kunz multiplicity and transcendental $F$-signature.
\end{mainresult}

The two multiplicities also give accumulation points of the kind considered by \c{C}elikba\c{s}, Dao, Huneke and Zhang \cite[p.~151]{CDHZ}. Adjoining a $2^r$th root of the defining equation realizes the normalized Frobenius colengths as multiplicities of hypersurfaces one dimension higher. The resulting dyadic rationals approach $\eHK(R_1)$ from below in dimension five and $\eHK(R_2)$ from below in dimension eight. For $R_2$, the coefficientwise square used in the proof also gives the successive differences and their rate of decay (Corollary~\ref{cor:monotone}). Finally, the series detect infinitely many indecomposable summands among the Frobenius pushforwards of each $R_m$. For $R_2$ the two limiting invariants also determine the signatures measuring $F$-rationality and the positive-degree Frobenius Betti numbers of its residue field.

\section{Branch matching and comparison of colengths}\label{sec:cubic}

The colengths in Theorem~\ref{main} are the same for every geometrically integral plane cubic with an ordinary node in characteristic two. Over an algebraic closure choose coordinates with node $[0:0:1]$ and tangents $x=0$, $y=0$. After scaling, its equation is $xyz+f_3(x,y)$. A linear change in $z$ removes the $x^2y$ and $xy^2$ terms. The coefficients of $x^3,y^3$ are both nonzero by irreducibility, so rescaling the coordinates gives $h$. Invertible linear substitutions preserve $(x^q,y^q,z^q)$, and extension of the coefficient field preserves colengths.

Since $h$ is defined over $\mathbb F_2$, we perform the calculation over that field. Write $P=\mathbb F_2[x,y,z]$, $B=P/(h)$ and $\mathbb D=[0,1]\cap\mathbb Z[1/2]$. For a function $g$ on $\mathbb D$, put
\[
 (T_0g)(t)=g(t/2),\qquad (T_1g)(t)=g((1+t)/2),\qquad \epsilon(t)=t(1-t).
\]
Monsky's functions $\phi_m$ vanish at $0$ and $1$. With $c_m=4m+3$, their values on the two halves of $\mathbb D$ satisfy
\begin{equation}\label{rules}
\begin{array}{c|cc}
 &8T_0\phi_m&8T_1\phi_m\\ \hline
 m=0&\phi_1+3t&\phi_0+3(1-t)\\
 m>0\text{ even}&\phi_{m+1}+c_m t&\phi_{m-1}+\epsilon+c_m(1-t)\\
 m\text{ odd}&\phi_{m-1}+\epsilon+c_m t&\phi_{m+1}+c_m(1-t).
\end{array}
\end{equation}
The right-hand side at each interior dyadic point uses values with smaller denominator. The endpoint values thus determine all the functions.

\subsection{The two branches at the node}\label{sec:node}

The normalization $\nu:\PP^1\to C$ is
\[
 [s:t]\longmapsto[s^2t:st^2:s^3+t^3].
\]
The two points $s=0$ and $t=0$ map to the node $p=[0:0:1]$. A section of $\OO_C(d)$ is a binary form of degree $3d$ with equal values at these points in the frames $t^{3d}$ and $s^{3d}$. For $d\ge0$ the parametrization and the Hilbert function of $B$ identify this space with $B_d$.

Fix $1\le j\le q$ and let $X_j=V(h^j)$. Every section of $\OO_{X_j}(D)$ lifts to a plane form since $H^1(\PP^2,\OO(n))=0$ for every $n$. Put $a=x+y$, $b=z+a$ and $P_{\mathrm f}=[1:1:0]$. The line $b=0$ cuts out the flex $P_{\mathrm f}$ with multiplicity three. On the chart $b=1$ we have $h=a^3+xy$. In the completion at $p$ choose $\theta$ with $\theta(0)=0$ and $\theta^2+\theta=a$. Then
\[
 u=x+a\theta,\qquad v=x+a(1+\theta),\qquad a=u+v,\qquad h=uv .
\]
Write $S=\mathbb F_2[[u,v]]$ and $T_j=S/(u^j,v^j)$.

We separate the branches of $X_j$ at $p$. For any power of two $n\ge j$, put $I_j=(x^n,y^n):h^{n-j}$. Since $z$ is a nonzerodivisor modulo $(x^n,y^n)$, it remains one modulo $I_j$. Thus $I_j$ is saturated. Its projective support is $p$. In the completed local ring, $(x^n,y^n)=(u^n,v^n)$, so $I_jS=(u^j,v^j)$. Let $\Dcal_j$ be the endomorphism algebra of the ideal sheaf induced by $I_j$ on $X_j$. Away from $p$ it is $\OO_{X_j}$. In $S/(u^jv^j)$ the ideal is the direct sum $(u^j)\oplus(v^j)$, with no homomorphisms between the two summands. Thus the completed algebra is $S/(u^j)\oplus S/(v^j)$, and the conductor sequence is
\begin{equation}\label{conductor}
 0\longrightarrow\OO_{X_j}\longrightarrow\Dcal_j
   \xrightarrow{\delta}T_j\longrightarrow0 .
\end{equation}
We regard $T_j$ as a sheaf at $p$. In characteristic two, $\delta$ adds the two branch values.

On $X_j$, let $Z_q$ be the Cartier divisor at $P_{\mathrm f}$ with local equation $(a/x)^q$. The identity $b^q(xy)^q=a^{3q}$ on $X_j$ gives $\operatorname{div}(b^q)=3Z_q$. Write $\rho_j:Y_j=\Spec_{X_j}\Dcal_j\to X_j$ for the curve obtained by separating the branches. The $q$th-power map $F_q:X_j\to C$ lifts on $Y_j$ to $\widetilde F_q:Y_j\to\PP^1$. Indeed $s/(s+t)=x/a$ on $C$, and $(x/a)^q$ is $1+\theta^q$ on one branch and $\theta^q$ on the other, since $q\ge j$. Away from $p$ the lift is given by the inverse normalization. As $(s+t)/s=a/x$ at the flex, $\widetilde F_q^*\OO_{\PP^1}(1)=\rho_j^*\OO_{X_j}(Z_q)$. On $T_j$ we have $\theta^q=0$, so the two restrictions of the lift are the constant maps to the points over $p$. We compare their values using $b=1$ and the frame $s+t$.

\begin{lemma}\label{common}
For $m\in H^0(\Dcal_j(D))$, the value $\delta(m)\in T_j$ is the restriction of a section of $\OO_{X_j}(D)$.
\end{lemma}
\begin{proof}
First we show that $h^j\in I_j^2$. For $j=2l+1$, invert $b$, adjoin $\theta$ with $\theta^2+\theta=a/b$, and put $w=x+a\theta$. Its norm is $h/b$. Write $b^l w^{2l+1}=G_l+K_l\theta$. The coefficients are plane forms. Indeed $(G_0,K_0)=(x,a)$, and multiplication by $bw^2$ acts by
\[
 \begin{pmatrix}bx^2+a^3&a^3\\ba^2&by^2+a^3\end{pmatrix}.
\]
At $p$ the element and its conjugate lie in $(u^j,v^j)$. Their sum is $K_l$, so both coefficients lie in $I_j$ by saturation. Taking norms gives $h^j=bG_l^2+bG_lK_l+aK_l^2$. For even $j$, square the corresponding expression for $h^{j/2}$. The local description gives $f^2\in I_j$ whenever $f\in I_{j/2}$.

Write $h^j=\sum_i f_i g_i$ with $f_i,g_i\in I_j$ homogeneous. By \eqref{conductor}, $f_i m$ and $g_i m$ descend to $X_j$. Choose plane lifts $F_i,G_i$. The two terms in $f_iG_i+g_iF_i$ cancel on $X_j$, so
\[
 L=\frac{\sum_i (f_iG_i+g_iF_i)}{h^j}
\]
is a plane form of degree $D$. To compute its value, set $U=u^j,V=v^j$ and write $f_i=c_iU+d_iV$, $g_i=e_iU+k_iV$. For representatives $m_s,m_t$ of the branch values we have $F_i=c_iUm_s+d_iVm_t$ modulo $UV$, and likewise for $G_i$. In $T_j$ this gives
\[
 L=(m_s+m_t)\sum_i(c_i k_i+d_i e_i).
\]
The sum is $1$ in $T_j$. Indeed it is the coefficient of $UV$ in $\sum_i f_i g_i=UV$ modulo $(U,V)^3$, and $U,V$ form a regular sequence. Thus $L|_{T_j}=\delta(m)$.
\end{proof}

\begin{lemma}\label{interchange}
Let $(\alpha,\beta)=(1,2)$ or $(2,4)$. If a section of $\Dcal_j(D)(-\beta Z_q)$ has branch values $(A,B)$ in $T_j$, a section of $\Dcal_j(D)(-\alpha Z_q)$ has values $(B,A)$.
\end{lemma}
\begin{proof}
For $(2,4)$, divide the given section by $b^q$ and apply Lemma~\ref{common}. Multiply its common lift by $b^q$ and add it to the original section. The two values are interchanged, and the added section vanishes along $3Z_q$. For $(1,2)$ and $j\le q/2$, use $b^{q/2}$ and $Z_q=2Z_{q/2}$ in the same argument.

For $(1,2)$ define
\[
 V_D=\im\!\left(H^0(\OO_{X_j}(D)(-Z_q))\longrightarrow T_j\right),
 \qquad W_D=\delta H^0(\Dcal_j(D)(-2Z_q)).
\]
It suffices to prove $W_D\subseteq V_D$. A common lift of $A+B$ changes $(A,B)$ to $(B,A)$ by addition. Set $D'=3j-3+q-D$. With $\lambda_j(r)=[u^{j-1}v^{j-1}]r$, multiplication followed by $\lambda_j$ is a perfect pairing on $T_j$. We claim that
\begin{equation}\label{duality}
 W_D=V_{D'}^\perp,\qquad \langle r,s\rangle=\lambda_j(rs).
\end{equation}
Twist \eqref{conductor} by $M=\OO_{X_j}(D)(-2Z_q)$. The boundary $\partial:T_j\to H^1(X_j,M)$ has kernel $W_D$. Since $X_j$ is a projective plane hypersurface, it is Cohen--Macaulay of pure dimension one. By Serre duality \cite{StacksCurve}, the source of the transpose is $H^0(X_j,\omega_{X_j}M^{-1})$. In the local frames fixed above, the transpose sends $s\in H^0(\omega_{X_j}M^{-1})$ to the functional $r\mapsto\lambda_j(r\,s|_{T_j})$. Here $\omega_{X_j}M^{-1}=\OO_{X_j}(D')(-Z_q)$, since $\omega_{X_j}=\OO_{X_j}(3j-3)$ and $\operatorname{div}(b^q)=3Z_q$.

To verify the functional in these frames, splice the completed conductor sequence with the resolution of $S/(UV)$ for $U=u^j,V=v^j$. This gives
\[
 0\longrightarrow S\xrightarrow{UV} S\xrightarrow{\rm diag}
 S/(U)\oplus S/(V)\xrightarrow{+}T_j\longrightarrow0 .
\]
A comparison from the Koszul resolution of $(U,V)$ sends $1$ to $(1,0)$ in degree zero and $e_1,e_2$ to $0,V$ in degree one. Since $d_2(1)=Ve_1+Ue_2$, it sends $1$ to $1$ in degree two. The extension class is therefore $1\in\operatorname{Ext}^2_S(T_j,S)=T_j$. The coordinate change preserves the volume form $da\wedge dx=du\wedge dv$, so its residue functional is $\lambda_j$. This identifies $\im(\partial^\vee)$ with $V_{D'}$ and proves \eqref{duality}. Interchanging $D,D'$ in \eqref{duality} gives $V_D=W_{D'}^\perp$. It remains to prove $\lambda_j(W_DW_{D'})=0$.

Each section used to define $W_D$ can be written $m_g=b^qg/a^q$, where
\[
 g\in S(q,j)_D,\qquad S(q,j)=(x^q,y^q,h^j):z^\infty .
\]
Indeed multiplying the section by $a^q/b^q$ gives a section regular at the flex. It descends at the node since $a^q$ kills $T_j$, and lifts to a plane form $g$. Conversely $m_g$ is regular on both branches precisely when, locally at $p$,
\[
 g\in(u^q,v^j)\cap(u^j,v^q)=(u^q,v^q,h^j).
\]
The saturation by $z$ imposes this condition at $p$. Write $A_g,B_g$ for the branch values of $m_g$. For $g,g'$ of degrees $D,D'$, the product $m_gm_{g'}/b^q$ is a separated canonical section, since the numerator vanishes along $4Z_q$. In the conductor sequence twisted by $\omega_{X_j}$, its difference of branch values has zero boundary. Pairing that boundary with $1\in H^0(\OO_{X_j})$ gives zero. Hence
\begin{equation}\label{diagonal}
 \lambda_j(A_gA_{g'}+B_gB_{g'})=0 .
\end{equation}

We prove the required vanishing by induction on $q$. For $j=q$, $S(q,q)=(x^q,y^q)$. For either generator times $c$, the difference of branch values is $b^q c$ on $T_j$, since the two values of $(x/a)^q$, and also those of $(y/a)^q$, differ by one. This gives a common lift vanishing along $3Z_q$.

For the induction step put $Q=2q$, $J=q+j$. Modulo $(x^{2q},y^{2q})$ we have $h^{q+j}=(xyz)^q h^j$. Multiplication by $(xy)^q$ identifies the monomials with $x,y$ exponents below $q$ with a direct summand over $\mathbb F_2[z]$ in the larger quotient. Saturation gives
\[
 S(Q,J)=(x^Q,y^Q)+(xy)^q S(q,j).
\]
A multiple of $x^Q$ or $y^Q$ gives a common lift as above and thus pairs to zero with the complementary $W$ by \eqref{duality}. It remains to pair $G=(xy)^qg$ with $G'=(xy)^qg'$.

On the branch $v^J=0$, use $b=1$, $a^{2q}=u^{2q}$ and $(xy)^q=(uv)^q+a^{3q}$. Dividing $G$ by $a^{2q}$ gives $a^qg+v^q\widehat A_g$, where $\widehat A_g$ represents $g/u^q$ modulo $v^j$. On the other branch the result is $a^qg+u^q\widehat B_g$. The common term cancels in $\delta m_G$. In the product with $\delta m_{G'}$, terms divisible by $u^{2q}$ or $v^{2q}$ vanish in $T_J$. The two cross terms give
\[
 \lambda_J(\delta m_G\,\delta m_{G'})
 =\lambda_j(A_gB_{g'}+B_gA_{g'})
 =\lambda_j\big((A_g+B_g)(A_{g'}+B_{g'})\big)
\]
by \eqref{diagonal}. The smaller degrees are complementary since $(3J-3+Q)-4q=3j-3+q$. The last expression vanishes by induction. Every pair with $Q/2<J<Q$ is obtained this way. The cases $J\le Q/2$ and $J=Q$ were treated above.
\end{proof}

\subsection{Comparison of two colength functions}\label{sec:compare}

For a homogeneous $\m$-primary ideal $I\subset P$, where $\m=(x,y,z)$, define
\[
 \Phi_I(j/q)=q^{-3}\dim_{\mathbb F_2}P/(I^{[q]},h^j).
\]
Flatness of Frobenius makes $\Phi_I$ well defined on $\mathbb D$. Choose homogeneous forms $f_i\in I$, of degrees $d_i$, which generate $IB$ as a graded ideal. Their syzygy bundle $S_I$ on $C$ is the kernel in
\begin{equation}\label{rowseq}
 0\longrightarrow S_I\longrightarrow\bigoplus_i\OO_C(-d_i)
     \xrightarrow{(f_i)}\OO_C\longrightarrow0 .
\end{equation}
Adding a redundant form of degree $d$ adds a summand $\OO_C(-d)$ to $S_I$. Graded generation also ensures that the $f_i^q$ generate $I^{[q]}$ modulo $h^j$ for $j\le q$.

Let $\Ccal_L(\alpha,\beta)$ be the bundle whose pullback to $\PP^1$ has $L$ copies of $\OO(-\alpha)$ followed by $L$ copies of $\OO(-\beta)$. In the frames $s$ and $t$, glue the fiber at $t=0$ of each summand to the fiber at $s=0$ of the next by the identity, including the last and first summands. We call this a cycle of half-length $L$, and put $\Bcal_L=\Ccal_L(1,2)$. A plane twist by $-N$ adds $3N$ to each index. Frobenius pullback multiplies the indices by $q$ and preserves the gluing \cite{Burban}.

For $1\le j\le q$, consider the pullback of a cycle $\Ccal$ by $F_q:X_j\to C$. Write $\delta_1,\ldots,\delta_{2L}$ for the indices of its summands in cyclic order. In the frames of Section~\ref{sec:node}, let $A_i$ and $B_i$ be the values of a section of the $i$th summand at the branches over $s=0$ and $t=0$. There is an exact sequence of sheaves on $X_j$
\begin{equation}\label{eq:matching}
\begin{split}
0\longrightarrow (F_q^*\Ccal)(D)
 &\longrightarrow \bigoplus_{i=1}^{2L}\Dcal_j(D)(-\delta_i Z_q)
 \xrightarrow{\ \Delta\ } T_j^{\,2L}\longrightarrow0,\\[-1pt]
 &\Delta((A_i,B_i)_i)=(B_i-A_{i+1})_i,\qquad A_{2L+1}=A_1 .
\end{split}
\end{equation}
Indeed, tensor the conductor sequence \eqref{conductor} with the locally free sheaf $(F_q^*\Ccal)(D)$.
Since $F_q\rho_j=\nu\widetilde F_q$, the middle term is the pushforward under $\rho_j$ of the normalized line summands pulled back to $Y_j$. This gives the displayed twists by $-\delta_i Z_q$. On the two copies of $\operatorname{Spec}T_j$ the lift is constant, with values $[0:1]$ and $[1:0]$. Thus the gluing of $\Ccal$ identifies the second value of summand $i$ with the first value of summand $i+1$ over the entire intersection. This gives $\Delta$. By left exactness of $H^0$, the sections of $(F_q^*\Ccal)(D)$ are precisely the tuples satisfying these matching conditions in degree $D$.

\begin{proposition}\label{compare}
Suppose $S_I$ and $S_J$ each have one cycle, with all other summands twists of $\OO_C$. For $r\in\mathbb Z$ and $c\in\mathbb Z_{\ge0}$, let their cycles be $\Ccal_L(\alpha,\beta)$ and $\Ccal_L(\alpha+3c,\beta+3c+3)$. If $(\alpha,\beta)=(3r+1,3r+2)$ or $(3r+2,3r+4)$, then
\[
 \Phi_J(t)-\Phi_I(t)=\bigl(\Phi_J(1)-\Phi_I(1)\bigr)t+\epsilon(t)
 \qquad(t\in\mathbb D).
\]
\end{proposition}
\begin{proof}
Fix $t=j/q$ with $1\le j\le q$ and put $R_j=P/(h^j)$. The pullback of \eqref{rowseq} to $X_j$ is exact because its quotient is locally free. Since $H^0(X_j,\OO_{X_j}(D))=(R_j)_D$, the row degrees and the sections of its kernel determine the Hilbert series of the quotient.

Consider the path of $L$ summands of index $\beta$ in the first cycle. Give it free exterior ends and impose the $L-1$ conditions $B_i=A_{i+1}$ along the path. Every such tuple is the restriction of a global section of the cycle. Indeed, reverse the order of its sections. For each of them Lemma~\ref{interchange} supplies a section on the opposite path with interchanged values. All $2L$ equalities in \eqref{eq:matching} now hold, including the two at the ends of the original path. The common shift $3r$ is absorbed in the plane twist by $\operatorname{div}(b^q)=3Z_q$. The same proof applies to the second cycle. After removing $3c$ from both indices, a section for $\beta+3$ is also a section for $\beta$ and has the same branch values. In either cycle the kernel of restriction is the low-index path with zero exterior values.

Let $K(T)$ be the Hilbert series of the first high-index path with free ends. In sufficiently large degree $D$, a section of each separated summand can take arbitrary values on the two copies of $\operatorname{Spec}T_j$. This follows by applying Serre vanishing to the ideal sheaf of their union on $Y_j$. The $L-1$ matching conditions are then independent. By \eqref{conductor}, separation adds $j^2$ to the Hilbert function of each summand in these degrees. The divisor $Z_q$ has length $qj$, since $3Z_q=\operatorname{div}(b^q)$ has length $3qj$. The twist $-\beta Z_q$ therefore changes the constant term of that Hilbert function by $-\beta qj$. Consequently
\[
 K(T)=L\Hilb_{R_j}(T)+\frac{j^2-L\beta qj}{1-T}+O(1)
 \qquad (T\to1).
\]

For the row of $I$ put $A_I(z)=1-\sum_i z^{d_i}+\sum_l z^{e_l}$, where $\OO_C(-e_l)$ are the line summands outside the cycle. Define $A_J$ similarly and set $z=T^q$. The low-index path for $J$ is the one for $I$ shifted by $cq$, and the high-index path is shifted by $(c+1)q$. Multiplication by the appropriate power of $b^q$ realizes both identifications and preserves branch values in the chart $b=1$. Subtract the quotient series for $I$, multiplied by $T^{cq}$, from that for $J$. The low-index paths cancel. The difference is
\[
 (A_J(z)-z^cA_I(z))\Hilb_{R_j}(T)+(z^{c+1}-z^c)K(T).
\]
Put $G(z)=A_J(z)-z^cA_I(z)+L(z^{c+1}-z^c)$. Rank and degree in \eqref{rowseq} give $A_I(1)=-2L$ and $3A_I'(1)=-L(\alpha+\beta)$. The corresponding identities for $J$ give $G(1)=G'(1)=0$. Since $\Hilb_{R_j}(T)=(1-T^{3j})/(1-T)^3$ and $z^{c+1}-z^c=-q(1-T)+O((1-T)^2)$, evaluation at $T=1$ gives the difference of colengths
\[
 \left(\tfrac32G''(1)+L\beta\right)q^2j-qj^2.
\]
After division by $q^3$ this is $\gamma t-t^2$. Its value at $t=1$ gives $\gamma=\Phi_J(1)-\Phi_I(1)+1$. At $t=0$ both colengths are zero.
\end{proof}

For example, take $I=\m$ and $J=\m^{[4]}$. The row $(x,y,z)$ has cycle indices $(4,5)$, as seen from its two columns in Section~\ref{sec:proof}. The fourth-power row has indices $(16,20)$. Since $E_1(1)=1$ and $E_4(1)=35$ by the formula for the reduced cubic \cite[Theorem~3]{BC}, the proposition gives
\begin{equation}\label{eq:four-example}
 E_{4q}(j)-E_q(j)=35q^2j-qj^2,\qquad 0\le j\le q .
\end{equation}
Here the term $-qj^2$ records the length of the intersection of the two branches.

\section{Colon ideals and syzygy cycles}\label{sec:colon}

The comparison in Proposition~\ref{compare} requires the cycle of each syzygy bundle. We determine the cycles for $(I^{[2]}:h)$ and $(I^{[4]}:h)$ in Theorem~\ref{transitions}.

\subsection{Extension classes on the double cubic}

For $c=2,4$ put $U_c(I)=(I^{[c]}:h)B/I^{[c]}B$. This graded module measures the forms that must be added to the row $(f_i^c)$. Use the $c$th-power map $F_c:X_2\to C$ and put $\widetilde V=F_c^*S_I$, $V=\widetilde V|_C$.

\begin{lemma}\label{transfer}
The row $(f_i^c)$ induces an injection $B_d/(I^{[c]}B)_d\hookrightarrow H^1(C,V(d))$. Under this injection $U_c(I)_d$ is the image of the boundary in
\begin{equation}\label{obseq}
 0\longrightarrow V(d)\xrightarrow{h}\widetilde V(d+3)
   \longrightarrow V(d+3)\longrightarrow0 .
\end{equation}
If forms $g_\nu$ of degrees $d_\nu$ represent such classes, the kernel $S'$ of the enlarged row is an extension
\[
 0\longrightarrow V\longrightarrow S'
   \longrightarrow\bigoplus_\nu\OO_C(-d_\nu)\longrightarrow0
\]
with extension classes the images of the $g_\nu$ in $H^1(C,V(d_\nu))$.
\end{lemma}
\begin{proof}
Lift the coefficients of a section of $V(d+3)$ to homogeneous polynomials $\widehat u_i$. Their relation on $C$ gives $\sum_i\widehat u_i f_i^c=hg$ for a form $g$ of degree $d$. Changing the lifts changes $g$ by an element of $(f_i^c)$. Conversely the $f_i^c$ generate $I^{[c]}$ modulo $h^2$, since $c\ge2$ and the $f_i$ generate $IB$. The condition $hg\in I^{[c]}+(h^2)$ depends only on $\bar g\in B$ and is equivalent to $\bar g\in(I^{[c]}:h)B$. Indeed, if $hg=f+h^2v$ with $f\in I^{[c]}$, then $h(g-hv)=f$. This identifies the possible classes of $g$ with $U_c(I)_d$.

Locally on $C$ write $g=\sum_i b_i f_i^c$. Lift $b_i$ to $X_2$ and replace $\widehat u_i$ by $\widehat u_i-hb_i$. These local lifts of the original syzygy differ on overlaps by the cocycle that represents the boundary of $g$ for the old row. This gives the identification in \eqref{obseq}. For an enlarged row, projection onto the new coefficient summands is surjective because the old row already generates $\OO_C$. Its local splittings are given by the same $b_i$ and hence have the stated extension classes.
\end{proof}

We record how this construction depends on the input bundle. For a bundle $S$ on $C$ put $V_S=(F_c^*S)|_C$ and write $\partial_{S,d}$ for the boundary in \eqref{obseq}. Suppose $S_1=S(-a)\oplus\mathcal L$, where $\mathcal L$ is a sum of plane twists. The inclusion of the first summand gives the commutative square
\begin{equation}\label{eq:transport}
\begin{array}{ccc}
 H^0(C,V_S(d+3)) &\xrightarrow{\ \partial_{S,d}\ }& H^1(C,V_S(d))\\
 \big\downarrow &&\big\downarrow\\
 H^0(C,V_{S_1}(d+ca+3)) &\xrightarrow{\ \partial_{S_1,d+ca}\ }&
 H^1(C,V_{S_1}(d+ca)).
\end{array}
\end{equation}
The other summand in the lower row has zero boundary. Indeed, every section of a plane twist on $C$ lifts to a plane form and hence to $X_2$. Thus the image below is the image above in its first summand. Lemma~\ref{transfer} realizes each transported class by a form of the actual colon. The classes can be adjoined at once since
\[
 \operatorname{Ext}^1_C\!\left(\bigoplus_\nu\OO_C(-d_\nu),V\right)
       =\bigoplus_\nu H^1(C,V(d_\nu)).
\]
The tuple determines the enlarged kernel up to equivalence of extensions. Its component in the pullback of $\mathcal L$ is zero, so this summand splits off.

We next compute the dimension of the full colon quotient in each degree. Write $b(s)=\dim B_s$, equal to $0,1,3s$ for $s<0,s=0,s>0$, respectively. Let $\sigma_L(v,w)$ be the dimension of the sections of a cycle having $L$ summands of degree $v$ followed by $L$ of degree $w$ on the normalization. For the nonzero degrees $v>w$ used below, matching endpoints gives
\begin{equation}\label{sections}
 \sigma_L(v,w)=
 \begin{cases}
 L(v+w),&v,w>0,\\
 Lv-1,&v>0>w,\\
 0,&v,w<0 .
 \end{cases}
\end{equation}
In the middle case the $L$ positive summands form a path with zero values at both ends. Put $h_{c,L}(s)=\sigma_L(3s-c,3s-2c)$ for $c=1,2,4$.

\begin{lemma}\label{profile}
Suppose $S_I$ is the direct sum of $\Bcal_L(-N)$ and twists of $\OO_C$, with $L\ge1$. The entries below give $\dim U_c(I)_{cN+s}$ at the four successive values of $s$. All other graded pieces vanish.
\begin{equation}\label{profiles}
\begin{array}{c|c|cccc}
 c&s&\multicolumn{4}{c}{\dim U_c(I)_{cN+s}}\\ \hline
 2&-2,-1,0,1&L-1&3L&3L&L-1\\
 4&-1,0,1,2&L&3L+3&3L+3&L
\end{array}
\end{equation}
\end{lemma}
\begin{proof}
Set $R'=P/I^{[c]}$. Multiplication by $h$ gives
\[
 0\longrightarrow U_c(I)(-3)\longrightarrow (R'/hR')(-3)
   \xrightarrow{h}R'/h^2R'\longrightarrow R'/hR'\longrightarrow0 .
\]
The quotient by $h^2$ is the Frobenius pullback of $P/(I^{[c/2]},h)$. Its Hilbert series is $(1+T)^3\Hilb_{P/(I^{[c/2]},h)}(T^2)$. In the row formula \eqref{rowseq}, the generator and line-summand terms cancel by $(1+T^3)\Hilb_B(T)=(1+T)^3\Hilb_B(T^2)$. The cycle gives
\begin{equation}\label{profilesum}
 \dim U_c(I)_{cN+s}=h_{c,L}(s+3)+h_{c,L}(s)
 -\sum_{\substack{0\le i\le3\\s+3-i\ {\mathrm{even}}}}
     \binom3i h_{c/2,L}\big((s+3-i)/2\big).
\end{equation}
Substituting \eqref{sections} gives the four entries of each row of \eqref{profiles}. Above either range all cycle degrees are positive. The linear terms cancel because $\sum\binom3i=4$ and $\sum i\binom3i=6$ on each parity class. Below the ranges every term is zero.
\end{proof}

\subsection{Syzygies of the model row}\label{sec:model-row}

We compute the two colon rows for a bundle with one cycle. The transport of extension classes in \eqref{eq:transport} will then give the result for an arbitrary input bundle of the same type. Fix $L\ge3$, set $d=3L$ and write $w=s/t$ on the normalization. In $\mathbb Z/d\mathbb Z$ put
\[
 \begin{gathered}\Lambda=\{1,\ldots,L\},\qquad p_j=L+1+2j\pmod d\ (0\le j\le L),\\
 \Pi=\{p_0,\ldots,p_L\},\qquad \Sigma=\Lambda\cup\Pi .\end{gathered}
\]
Thus $p_L=1$. Our model row consists of $1+w^d$ and the monomials $w^p$ for $p\in\Sigma$, all of plane degree $L$. Let $I\subset P$ be generated by homogeneous lifts of the row and by $h$. The row has no common zero on $C$, so $I$ is $\m$-primary.

Let $K$ be the kernel of a surjective row $\bigoplus_{\ell=1}^{r}\OO_C(-e_\ell)\to\OO_C$. A homogeneous syzygy of binary degree $n$ is a column whose $\ell$th entry has degree $n-3e_\ell$. It gives a map $\OO_{\mathbb P^1}(-n)\to\nu^*K$. We write $v\to v'$ for its values at $s=0$ and $t=0$. We identify both fibers with $K_p$ using the plane frames from $b$. Indeed $\nu^*b=(s+t)^3$, so these are the frames used in the column calculations.

\begin{basislemma}
Suppose $r-1$ columns of degrees $n_i$ have endpoint pairs $v_i\to v_{\sigma(i)}$, where $(v_i)$ is a basis of $K_p$ and $\sigma$ is a permutation. If $\sum_i n_i=3\sum_\ell e_\ell$, then $K$ is the direct sum of the bundles associated to the orbits of $\sigma$. An orbit has the normalized summands $\OO_{\mathbb P^1}(-n_i)$ in their cyclic order. In particular, a one-element orbit with $n_i=3e$ gives $\OO_C(-e)$.
\end{basislemma}
\begin{proof}
The columns define a map $\bigoplus_i\OO_{\mathbb P^1}(-n_i)\to\nu^*K$ of equal-rank bundles. Its determinant is nonzero at $s=0$ and has degree zero by the degree sum. Hence it is nowhere zero, and the map is an isomorphism.

For an orbit $\Gamma$, let $W_\Gamma$ be the span of the $v_i$ with $i\in\Gamma$. The corresponding sum of normalized line bundles has fiber $W_\Gamma$ at both points over $p$. Since $K_p=\bigoplus_\Gamma W_\Gamma$, equality of the two branch values is equivalent to equality in each $W_\Gamma$. The conductor sequence for $j=1$, tensored with $K$, therefore splits the sheaf of matching sections into the direct sum over $\Gamma$. Within an orbit the endpoint pairs give the asserted cyclic gluing. On a one-element orbit of degree $3e$ the two frames agree, which is the gluing of the plane twist $\OO_C(-e)$.
\end{proof}

A column of degree $n'<n$ can be added to one of degree $n$ after multiplication by a binary form of degree $n-n'$. The monomials $t^{n-n'}$ and $s^{n-n'}$ prescribe the correction independently at the two endpoints.

In the model, consecutive monomials have gaps $g=1$ or $2$. The relation at a gap from $p$ to $p+g$ has coefficients $s^g,t^g$ and endpoint pair $e_{p+g}\to e_p$, where $e_p$ is the monomial coordinate. At the gap from $d-1$ to $1$ the additional coefficient of $1+w^d$ is $st$, since $w(1+w^d)=w+w^{d+1}$. It vanishes at both endpoints. The $2L$ monomial coordinates form a basis of the fiber. The pairs run once around it, with $L$ consecutive gaps of each size. The column degrees sum to $(2L+1)d$, the sum of the row degrees. The lemma gives the model kernel $\Bcal_L(-L)$.

\begin{proposition}\label{prop:model-rows}
For $c=2,4$, one can adjoin forms from $(I^{[c]}:h)B$ to the row of old powers so that its kernel has one cycle and the fixed columns listed below. In each pair the first entry is the degree and the second is the multiplicity. Row entries have plane degrees and columns have binary degrees. Write $D=3cL$.
\[
\setlength{\arraycolsep}{3pt}
\begin{array}{c|c|c|c}
 c & \text{new row entries} & \text{cycle columns} & \text{fixed columns}\\ \hline
 2 & (2L-2,L-1),(2L-1,L)
   & (D-4,L-2),(D+1,L-2) & (D,2L+3)\\
 4 & (4L-1,L),(4L,L+1)
   & (D+1,L+2),(D+2,L+2) & (D,L-2),(D+6,L-1)
\end{array}
\]
In each cycle the columns of either degree occur consecutively.
\end{proposition}

To construct forms in either colon we use \eqref{obseq}. For a lifted relation $\sum\widehat u_i f_i^c=hg$, differentiation in $z$ gives
\begin{equation}\label{derivative}
 xy\,\bar g=\sum_i\overline{\partial_z\widehat u_i}\,\bar f_i^{\,c}\quad\text{in }B .
\end{equation}
The same criterion applies to both rows. Equation~\eqref{derivative} identifies the new entries from lifted relations. For the enlarged row we construct homogeneous columns, check their syzygies and exhibit a basis of the fiber at the node permuted by their endpoint values. Equality of the two degree sums then gives the kernel bundle by the basis lemma. Generation of the full colon ideal is checked separately in Theorem~\ref{transitions}.

\subsubsection{Fourth powers}\label{sec:cubic:four}
The old fourth powers have binary degree $D=12L$. Keep the forms $Z=1+w^D$ and $F_p=w^{4p}$ for $p\in\Sigma$. On a gap of size one, multiply the fourth power of its column by $st$. The two coefficients $s^5t,st^5$ lift to $xz+y^2,yz+x^2$. By \eqref{derivative} this gives
\[
 w^{4p-1}(1+w^3),\qquad p\in\Lambda ,
\]
in binary degree $D-3$.

There is one further form $M_j$ of degree $D$ at each $p=p_j$. Let $g_-,g_+$ be the gaps to its neighbours. Multiply the adjacent fourth-power columns by $t^{9-4g_-}$ and $s^{9-4g_+}$ and add them. At $p$ the coefficient is $1+w^9$, with neighbouring coefficients $w^{4g_-}$ and $w^{9-4g_+}$. On $C$ the required lifts and their $z$-derivatives are
\[
\begin{array}{c|ccccc}
 \text{coefficient}&1+w^9&w^4&w^5&w^8&w\\ \hline
 \text{lift}&z^3+x^3+y^3&xy^2&x^2y&xz^2+y^2z+x^2y&yz^2+x^2z+xy^2\\
 \partial_z&z^2&0&0&y^2&x^2 .
\end{array}
\]
At the gap across zero the extra coefficient belongs to the column of the old generator $Z$. Its derivative in \eqref{derivative} is therefore multiplied by $Z$ and vanishes modulo $Z$. Since $xy=w^3$ is a unit modulo $Z=1+w^D$ on $t=1$, we may compute \eqref{derivative} modulo $Z$. It gives
\begin{equation}\label{fourrow}
 M_j=w^{4p-3}+w^{4p+3}+
           \mathbf1_{j>0}w^{4p-9}+\mathbf1_{j<L}w^{4p+9}\pmod{Z}.
\end{equation}
Here $\mathbf1$ denotes an indicator, and exponents are taken modulo $D$. This changes the form by at most a scalar multiple of the retained generator $Z$.

To compute its syzygies, work in binary degree $D+k$. All row entries except $Z$ vanish at both endpoints. If the dehomogenized sum of the other terms vanishes modulo $Z=1+w^D$, its quotient by $Z$ has degree at most $k$. Homogenizing that quotient gives the missing coefficient of $Z$. Its two endpoint values vanish because the other row entries vanish there and $Z$ has value one in both frames. Work in $\mathbb F_2[Y]/(Y^d+1)$ with $Y=w^4$. For the degree-$(D-3)$ forms, the $F_p$, and the $M_j$, collect the coefficients of $w^i$ in polynomials $A_i,B_i,C_i$ supported on $\Lambda,\Sigma,\Pi$, respectively. The ranges of $i$ are $0,\ldots,k+3$ for $A_i$ and $0,\ldots,k$ for $B_i,C_i$. If $C|_p$ denotes its single term at $p$, set
\[
 D_1(C)=Y^{-1}C+Y^2(C-C|_{p_L}),\qquad
 D_3(C)=C+Y^{-3}(C-C|_{p_0}).
\]
By \eqref{fourrow} the condition for a column is
\begin{equation}\label{fourres}
 \sum_i w^i\big((w^3Y^{-1}+w^2)A_i+B_i+wD_1(C_i)+w^3D_3(C_i)\big)=0
 \pmod{1+w^D}.
\end{equation}

In degree $D$, let $A$ have support in $\{2,\ldots,L-1\}$. Take $A_0=A_3=A$, $A_1=Y^{-1}A$, $A_2=YA$, $B_0=(Y^{-1}+Y^2)A$ and $C_0=0$. These solve \eqref{fourres} and give $L-2$ fixed columns with endpoint vector $(A,(Y^{-1}+Y^2)A,0)$.

For degree $D+1$, set $C_0=C_1=0$. Take $A_0,A_4,T$ with supports in $[2,L]$, $[1,L-1]$, $[1,L]$, respectively, and put
\[
 A_1=Y^{-1}A_0+A_4,\quad A_2=T,\quad A_3=A_0+YA_4,\qquad
 B_0=A_1+YT,\quad B_1=T+YA_3 .
\]
Each of the four residue classes of $w$ in \eqref{fourres} vanishes. On endpoint vectors with $A$ supported on $\Lambda$ and $B$ on $[1,L+1]$, define $\pi(A,B)=B+(Y^{-1}+Y^2)A$. Its kernel is the span of the fixed columns above, and its image has basis $1,Y,\dots,Y^{L+2}$. After applying $\pi$, the right endpoint is $T+YA_0+Y^{-1}A_4$ and the left endpoint is $Y$ times it. The support of $Y^{-1}A_4$ includes exponent zero, $T$ supplies exponents $1,\ldots,L$, and $YA_0$ supplies exponent $L+1$. By choosing one monomial at a time, the right endpoint therefore runs through $1,Y,\ldots,Y^{L+1}$. The resulting columns give the path
\[
 Y^{L+2}\longrightarrow Y^{L+1}\longrightarrow\cdots\longrightarrow 1 .
\]
For $1$ choose the representative $(Y,Y^3)$, for $Y^{L+2}$ choose $(Y^L,Y^{L-1})$, and for $Y^i$ with $1\le i\le L+1$ choose $(0,Y^i)$. Any discrepancy between a column endpoint and its chosen representative lies in $\ker\pi$. A fixed column with endpoints $(k,k)$ gives $(k,0)$ after multiplication by $t$ and $(0,k)$ after multiplication by $s$. Since the fixed columns span $\ker\pi$ at each end, we may use the same representative $v_j$ on both edges meeting at $Y^j$.

\begin{correctionlemma}
Let $E$ be the endpoint space with $A$ supported on $[1,L]$ and $B$ on $[1,L+1]$. Write $v_j$ for the chosen lift of $Y^j$ and $K=\ker\pi$. In degree $D+2$, multiples of the fixed columns of degree $D$ and the path columns of degree $D+1$ give precisely the endpoint corrections
\[
 \bigl(K+\langle v_1,\ldots,v_{L+2}\rangle\bigr)
 \times
 \bigl(K+\langle v_0,\ldots,v_{L+1}\rangle\bigr)
 \quad\text{in }E\times E .
\]
\end{correctionlemma}
\begin{proof}
A path column has endpoints $(v_{j+1},v_j)$ for $0\le j\le L+1$. Multiplying by $t$ or $s$ gives $(v_{j+1},0)$ or $(0,v_j)$ in the endpoint frames. A fixed column has endpoints $(k,k)$; multiplication by $t^2$ or $s^2$ gives $(k,0)$ or $(0,k)$. These $k$ span $K$. This proves that every pair in the displayed product occurs. Conversely a homogeneous multiplier acts at each end by its value there, so every correction from these columns lies in the product.
\end{proof}

The return path has degree $D+2$. Set $C_1=0$. Write $\alpha_j=(C_0)_{p_j}$ and $\beta_j=(C_2)_{p_j}$. Choose all the $\alpha_j$ and $\beta_L$ freely and put $\beta_j=\alpha_{j+1}$ for $j<L$. Take $A_0=\beta_LY$, $A_5=\alpha_0Y^L$, $A_2=A_3=A_4=0$. Equation \eqref{fourres} is satisfied by
\begin{align*}
 A_1&=Y^{-1}A_0+YA_5+D_3(C_0)+D_1(C_2),\\
 B_0&=A_1+YA_5,\qquad B_1=D_1(C_0)+YD_3(C_2),\qquad B_2=A_0 .
\end{align*}
To check the supports, the indices outside $\Sigma$ are $m_j=p_j+1$ for $0\le j<L$, with $m_{L-1}=0$. If $C_{p_j}=c_j$, then $D_1(C)_{m_j}=c_{j+1}$ and $(YD_3(C))_{m_j}=c_j$. This gives $(B_1)_{m_j}=0$. In $D_3(C_0)+D_1(C_2)$ the coefficients at $p_j$, $1\le j<L$, are $\alpha_j+\beta_{j-1}$, and those at $m_j$, $j<L-1$, are $\alpha_{j+2}+\beta_{j+1}$. All vanish. The two remaining coefficients outside $\Lambda$, at $0$ and $L+1$, cancel with $Y^{-1}A_0$ and $YA_5$.

The endpoint correction lemma leaves only the constant coefficient of $\pi$ on the left and the coefficient of $Y^{L+2}$ on the right. The new endpoints have no $F_{p_j}$ coordinate for $1\le j<L$. In the remaining endpoint coordinate followed by the $C$ coordinates, the pair is
\[
 (\beta_L,\alpha_0,\ldots,\alpha_L)\longmapsto
               (\alpha_0,\alpha_1,\ldots,\alpha_L,\beta_L).
\]
Writing $c_j$ for the unit vector of $M_j$, we obtain the path
$1\to c_L\to c_{L-1}\to\cdots\to c_0\to Y^{L+2}$. It closes the cycle, with $L+2$ columns of each degree $D+1,D+2$.

Finally, for $p=p_j$, $1\le j<L$, equation \eqref{fourrow} gives
\[
 (1+w^6)F_p+w^3 M_j+w^2F_{p_{j-1}}+w^4F_{p_{j+1}}=0\pmod Z .
\]
These $L-1$ columns of degree $D+6$ have the same unit vector of $F_p$ at both ends. We now have a single basis at the node. On the coordinates $A$ and $F_1,\ldots,F_{L+1}$ take the $L-2$ vectors spanning $\ker\pi$ and the chosen lifts of $1,Y,\ldots,Y^{L+2}$. Their independence follows by first applying $\pi$. Adjoin the $L+1$ vectors $c_j$ and the $L-1$ interior vectors $e_{F_{p_j}}$, which occupy disjoint new coordinates. All row entries except $Z$ vanish at the node, so these $4L+1$ vectors form a basis of the kernel fiber. In this basis the two paths form one orbit. The $\ker\pi$ vectors and the interior $F_{p_j}$ vectors are fixed. The column degrees sum to $(4L+1)D+9L$, the sum of the binary row degrees since $D=12L$. The basis lemma gives the cycle and both groups of line summands in the fourth-power row.

\subsubsection{Squares}\label{sec:square}
For squares the old binary degree is $D=6L$. At an interior vertex $p$, $2\le p\le L$, take the two adjacent gap-one columns. Their squares give a matched section with plane coefficients $x,z,y$ at the positions $p-1,p,p+1$. By \eqref{derivative} its image is $w^{2p-3}$ in binary degree $D-6$. On a gap of size two from $p$ to $p+2$, multiply the squared column by $st$. As before, its coefficients lift to $xz+y^2,yz+x^2$. They give $w^{2p-1}(1+w^3)$ in degree $D-3$. At the gap from $3L-1$ to $1$ this gives $1+w^{D-3}$. The extra coefficient on the old binomial lifts to $xy$ and has zero $z$-derivative.

Put $n=L-1$ and $b_0=2L+1$. For $0\le i<n$ the row reduces to
\begin{equation}\label{squarerow}
\begin{array}{c|c}
 \text{form}&\text{binary degree}\\ \hline
 A_i=w^{1+2i}&D-6\\
 B_i=w^{b_0+4i}(1+w^3),\quad Z=1+w^{D-3}&D-3\\
 C_i=w^{b_0+4i+5}&D .
\end{array}
\end{equation}
Indeed the old squares $w^{2p}$ for $1\le p\le L+1$ are multiples of the $A_i$ by $w,w^3$ or $w^5$, all plane forms of degree two. The old interior squares along the other path are the $C_i$. To recover $1+w^D$, use
\[
 (1+w^3)Z=1+w^D+w^3+w^{D-3},\qquad
 wB_{n-1}=w^{D-6}+w^{D-3}.
\]
Here $w^3=w^2 A_0$ and $w^{D-6}$ is an old square. Thus restoring the $L+2$ omitted forms adds $L+2$ copies of $\OO_C(-2L)$.

We use $\mathbf a_i,\mathbf b_i,\mathbf c_i$ for the unit vectors of the coordinates $A_i,B_i,C_i$. In degree $D-4$, the relation $w^2A_i+A_{i+1}=0$ gives $\mathbf a_{i+1}\to\mathbf a_i$. These $n-1$ columns give the first half of the cycle. To find fixed columns and the return path we work modulo $Z$ in degree $D+k$, for $k=0,1$. All forms in \eqref{squarerow} except $Z$ vanish at the endpoints, so the coefficient of $Z$ again has zero endpoint values.

Set $Y=w^4$. Let $P_r(Y)$ be the polynomial of degree less than $n$ whose $i$th coefficient is the coefficient of $w^r$ multiplying $B_i$. Define $Q_r$ in the same way for $C_i$. The contribution of these two groups, divided by $w^{b_0}$, is
\begin{equation}\label{squarearray}
 S=(1+w^3)\sum_{r=0}^{3+k}w^r P_r(Y)+w^5\sum_{r=0}^k w^r Q_r(Y).
\end{equation}

\begin{liftcriterion}
Fix $k=0$ or $1$ and the coefficients of the $B_i,C_i$ in \eqref{squarearray}. They extend to a homogeneous syzygy of binary degree $D+k$ if and only if
\[
 [w^j]S=0\qquad(k+3\le j\le4n).
\]
The left $A_0$ and right $A_{n-1}$ coordinates are then forced to be
\[
 a=[w^{4n+1}]S,\qquad a'=[w^{k+2}]S .
\]
All remaining $A$ coordinates at the two endpoints may be prescribed independently.
\end{liftcriterion}
\begin{proof}
The coefficient of $A_i=w^{1+2i}$ has degree at most $6+k$ on $t=1$. These terms therefore span the monomials of exponents $1$ through $M=b_0+k+2$. Only $A_0$ reaches the first exponent and only $A_{n-1}$ reaches the last. Put $N=D-3=b_0+4n$. Modulo $Z=1+w^N$, reduction sends an exponent $N+j$ to $j$. Since $\deg S\le4n+2+k$ and $2+k<b_0$, the terms that wrap around have exponents at most $2+k$ after reduction. Thus the forbidden exponents of $w^{b_0}S$ are $M+1,\ldots,N$, with $N$ representing the constant term. These give exactly the stated vanishing conditions. The coefficients at exponents $1$ and $M$ give $a$ and $a'$.

Choose coefficients for the $A_i$ to cancel the reduced polynomial. The sum before reduction is divisible by $Z$ and has degree at most $D+k$. Its quotient by $Z$ has degree at most $3+k$, the allowed degree for the last coefficient. Homogenizing the entries gives the required syzygy. Finally, the column $w^2A_i+A_{i+1}=0$ has endpoint pair $\mathbf a_{i+1}\to\mathbf a_i$. Multiplying it by $t^{4+k}$ or $s^{4+k}$ changes only the left or the right endpoint, respectively. For $0\le i<n-1$ these operations account for every free coordinate.
\end{proof}

For fixed columns take $k=0$. Write $[\,\cdot\,]_{<n}$ for truncation at degree $n$ and $\ell P=(P-P(0))/Y$. Given $P$ of degree less than $n$ and $a,a'\in\mathbb F_2$, set
\[
 \begin{gathered}
 P_0=P_3=P,\qquad P_1=\ell P,\qquad P_2=a'+[YP]_{<n},\\
 Q_0=\ell^2P+[YP]_{<n}+a'+aY^{n-1}.
 \end{gathered}
\]
The four residues of \eqref{squarearray} are $P+YP_1$, $P_1+Y(P_2+Q_0)$, $P_2+YP_3$, and $P_3+P$. They satisfy the cancellation condition and give the prescribed $a,a'$. The $B$ and $C$ coordinates are equal at both ends. Choose the same $A$ vector $a\mathbf a_0+a'\mathbf a_{n-1}$ at both ends as well. For $P=0$ the two unit choices of $a,a'$ give the fixed vectors
\[
 \mathbf a_0+\mathbf c_{n-1},\qquad \mathbf a_{n-1}+\mathbf c_0 .
\]
For $P=Y^i$ and $a=a'=0$ we obtain the fixed vector
$\mathbf b_i+\mathbf c_{i-2}+\mathbf c_{i+1}$, omitting indices outside $[0,n-1]$. This gives $n+2$ fixed columns in degree $D$.

For the return path take $k=1$. It suffices to choose $Q_1$ freely and put
\[
 P_0=P_1=P_4=0,\qquad P_3=a',\qquad
 P_2=[YQ_1]_{<n}+a'Y,\qquad Q_0=P_2+aY^{n-1}.
\]
The four residues of $S$ are now $0$, $aY^n$, $P_2+Y(a'+Q_1)$, and $a'$. The third is a multiple of $Y^n$, so the cancellation condition holds. The two $C$ vectors are $Q_0,Q_1$, and both $B$ vectors vanish. For $n>2$ the following choices give the return path.
\[
\begin{array}{ccc|c}
 Q_1&a&a'&\text{endpoints}\\ \hline
 Y^{n-2}&1&0&\mathbf a_0\to\mathbf c_{n-2}\\
 Y^{i-1}&0&0&\mathbf c_i\to\mathbf c_{i-1}\quad(2\le i\le n-2)\\
 0&0&1&\mathbf c_1\to\mathbf a_{n-1}.
\end{array}
\]
In the first row $Q_0=0$ because the two terms $Y^{n-1}$ cancel. In the middle row $Q_0=Y^i$, and in the last row $Q_0=Y$. The lifting lemma gives the indicated $A$ coordinates in each case. For $L=3$ we have $n=2$. Take $Q_1=0$ and $a=a'=1$. Then $Q_0=0$, giving the single edge $\mathbf a_0\to\mathbf a_1$. In every case the cycle has $n-1$ columns of degree $D-4$ and $n-1$ of degree $D+1$.

The vertices of this cycle are $\mathbf a_0,\ldots,\mathbf a_{n-1}$ and $\mathbf c_1,\ldots,\mathbf c_{n-2}$. Together with the $n+2$ fixed vectors displayed above they form a basis of the kernel fiber. Indeed, in a relation among them the $\mathbf b_i$ coordinates first eliminate $n$ fixed vectors. The coordinates $\mathbf c_0,\mathbf c_{n-1}$ then eliminate the other two. The cycle vertices are distinct coordinate vectors, also when $n=2$. The resulting $3n$ vectors exhaust the fiber. Their endpoint permutation has the stated cycle and $n+2$ one-element orbits. The degrees of the columns sum to $(n+2)D+(n-1)(2D-3)$, the sum of the binary row degrees in \eqref{squarerow}. The basis lemma therefore splits the reduced kernel over $C$. Restoring the $L+2$ redundant entries adds $L+2$ plane twists and gives $2L+3$ summands $\OO_C(-2L)$ in the full square row. This proves Proposition~\ref{prop:model-rows}.

\subsection{The change of cycle under Frobenius}\label{sec:model}

Proposition~\ref{prop:model-rows} determines the kernel for the model row. We now transfer its extension classes to the row of $I$ and compare Hilbert functions to prove that the resulting generators give the full colon ideal on $C$.

\begin{theorem}\label{transitions}
Let $S_I=\Bcal_L(-N)\oplus\mathcal L$, where $L\ge3$ and $\mathcal L$ is a sum of twists of $\OO_C$. For $c=2,4$ there is a homogeneous generating row of $(I^{[c]}:h)B$ containing the $c$th powers of the original row. If $S_c$ denotes its syzygy bundle and $\mathcal L^{[c]}$ the Frobenius pullback, then
\begin{align*}
 S_4&\simeq \Bcal_{L+2}(-4N)\oplus\OO_C(-4N)^{L-2}
               \oplus\OO_C(-4N-2)^{L-1}\oplus\mathcal L^{[4]},\\
 S_2&\simeq \Ccal_{L-2}(6N-4,6N+1)\oplus\OO_C(-2N)^{2L+3}\oplus\mathcal L^{[2]}.
\end{align*}
For $c=4$ the additional generators have degrees $4N-1$ and $4N$, with multiplicities $L$ and $L+1$. For $c=2$ their degrees are $2N-2$ and $2N-1$, with multiplicities $L-1$ and $L$.
\end{theorem}

\begin{proof}
Apply Proposition~\ref{prop:model-rows} to the model with kernel $\Bcal_L(-L)$. Its new row entries give a tuple of classes in the image of \eqref{obseq}. Apply \eqref{eq:transport} with $a=N-L$. The transported forms lie in $(I^{[c]}:h)B$, their degrees increase by $c(N-L)$, and the kernel of the enlarged row has the displayed decomposition.

It remains to prove that the transported forms generate the colon in every degree. If $K_c$ is the ideal of their row, including the old powers, then
\[
 I^{[c]}B\subseteq K_c\subseteq(I^{[c]}:h)B .
\]
Write $v_c(s)=\dim (K_c/I^{[c]}B)_{cN+s}$. Compute this difference from the two rows by \eqref{rowseq}. Their common line summands cancel. The additional generators and the change of cycle give
\begin{align*}
 v_4(s)&=L b(s+1)+3b(s)-(L-1)b(s-2)
           +h_{4,L}(s)-h_{1,L+2}(s),\\
 v_2(s)&=(L-1)b(s+2)+L b(s+1)-(2L+3)b(s)
           +h_{2,L}(s)-\sigma_{L-2}(3s+4,3s-1).
\end{align*}
For example, the coefficient $3$ of $b(s)$ in the first line counts $L+1$ new generators of degree $4N$ minus $L-2$ line summands of that degree. Substitution from \eqref{sections} gives
\begin{align*}
 \sum_s v_2(s)T^{s+2}&=(1+T)\bigl((L-1)(1+T^2)+(2L+1)T\bigr),\\
 \sum_s v_4(s)T^{s+1}&=(1+T)\bigl(L(1+T^2)+(2L+3)T\bigr).
\end{align*}
Their coefficients are precisely the dimensions of $U_2(I)$ and $U_4(I)$ in \eqref{profiles}, including at $L=3$. Thus the inclusion $K_c/I^{[c]}B\subseteq U_c(I)$ is an equality in each degree. This proves $K_c=(I^{[c]}:h)B$.
\end{proof}

\section{Monsky's formula and the signature of the pair}\label{sec:proof}

\subsection{The dyadic recursion}

Put $N_k=4^k$ and
\[
 A_k=(x^{N_k},y^{N_k},z^{N_k}):h^{(N_k-1)/3}.
\]
Flatness of Frobenius gives $A_{k+1}=(A_k^{[4]}:h)$. For $A_0=\m$, the row $(x,y,z)$ has normalized columns $(t,s,0)$ and $(s^2,t^2,st)$, of degrees $4$ and $5$. They form the cycle $\Bcal_1(-1)$.

\begin{lemma}\label{initial}
For every $k\ge0$, $A_kB$ has a generating row with syzygy bundle $\Bcal_{2k+1}(-N_k)$ up to direct summands that are twists of $\OO_C$.
\end{lemma}
\begin{proof}
For the first step put $I=(x^4,y^4,z^4)$ and $w=s/t$ on the normalization. The seven forms below lie in $I:h$, as multiplication by $h$ shows. The last column gives their restrictions on $t=1$.
\[
\begin{array}{c|l|l}
 i & g_i& g_i(w)\\ \hline
 0&x^3+y^3&w^3+w^6\\
 1&z^4&1+w^{12}\\
 2&y^4&w^4\\
 3&xy^3&w^5\\
 4&yz^3+xy^2z+x^2z^2&w+w^4\\
 5&y^2z^2+xz^3+x^2yz&w^8+w^{11}\\
 6&y^3z+x^2y^2&w^3
\end{array}
\]
Their restrictions generate an ideal containing $IB$, since $w^8=w^5+w^2(w^3+w^6)$. The columns below are syzygies after pullback to the normalization. Homogenize each entry to the column degree minus the binary degree of the corresponding form.
\[
\begin{array}{c|l}
 \text{degree}&\text{column}\\ \hline
 13&(1,0,0,w,0,0,1)\\
 13&(0,0,w,1,0,0,0)\\
 13&(0,0,1,0,0,0,w)\\
 14&(0,0,0,w,w^2,0,1)\\
 14&(w+w^4,w,0,0,1,w^2,0)\\
 14&(w^5,0,w,0,0,1,w^2)
\end{array}
\]
Write $e_i$ for the coordinate vector corresponding to $g_i$. The endpoint vectors run from $e_0+e_6$ through $e_3,e_2,e_6,e_4,e_5$ and back to $e_0+e_6$. The six vectors are independent. The degree sum $81=9+6\cdot12$ equals the sum of the row degrees, so the columns form a basis on $\PP^1$. Their endpoint permutation is the single cycle just displayed. By the basis lemma in Section~\ref{sec:model-row} they give $\Bcal_3(-4)$ on $C$. In degree $4+s$, the dimension of the ideal generated by this row modulo $IB$ is
\[
 b(s+1)+3b(s)+h_{4,1}(s)-h_{1,3}(s).
\]
By \eqref{sections} this is the $c=4$, $L=1$ row of \eqref{profiles}, so the forms generate $(I:h)B=A_1B$ in every degree. The fourth-power case of Theorem~\ref{transitions} now proves the lemma by induction.
\end{proof}

\begin{proof}[Proof of Theorem~\ref{main}]
For every homogeneous $\m$-primary ideal $I$ we have
\begin{equation}\label{halves}
 8T_0 \Phi_I= \Phi_{I^{[2]}},\qquad
 8T_1 \Phi_I= \Phi_{I^{[2]}}(1)+\Phi_{(I^{[2]}:h)} .
\end{equation}
The first equality follows by taking $t=j/q$. For the second, filter $P/(I^{[2q]},h^{q+j})$ by the image of $h^q$. The quotient by this image has length $q^3\Phi_{I^{[2]}}(1)$. The image is isomorphic, up to a degree shift, to $P/((I^{[2q]}:h^q),h^j)$. Flatness gives $(I^{[2q]}:h^q)=(I^{[2]}:h)^{[q]}$.

Write $P_k=\Phi_{A_k}$, $Q_k=\Phi_{A_k^{[2]}}$, and put $p_k=P_k(1)$, $r_k=Q_k(1)$ and $s_k=\Phi_{A_k^{[4]}}(1)$. The exact sequence in the proof of Lemma~\ref{profile} gives
\[
 \dim U_c(A_k)=2\Phi_{A_k^{[c]}}(1)-8\Phi_{A_k^{[c/2]}}(1).
\]
Since $L=2k+1$ for $A_k$, summing the two rows of \eqref{profiles} yields
\begin{equation}\label{slopes}
 r_k-4p_k=4L-1=8k+3,\qquad s_k-4r_k=4L+3=8k+7 .
\end{equation}
Subtract the linear interpolation between the endpoint values by setting
\[
 \psi_{2k}(t)=P_k(t)-p_kt,\qquad
 \psi_{2k+1}(t)=Q_k(t)-r_kt .
\]
These functions vanish at $0$ and $1$. We show they satisfy \eqref{rules}.

First, \eqref{halves} gives $8T_0P_k=Q_k$ and $8T_1Q_k=s_k+P_{k+1}$. Evaluation at $1$ in the second equality gives $p_{k+1}=8r_k-s_k$. It follows that
\[
 8T_0\psi_{2k}=\psi_{2k+1}+(8k+3)t,\qquad
 8T_1\psi_{2k+1}=\psi_{2k+2}+(8k+7)(1-t).
\]

For the two steps that lower the index, compare cycles with the same half-length. The indices for $A_k$ are $\alpha=3N_k+1$ and $\beta=3N_k+2$. In $A_k^{[4]}$ they increase by $3\alpha$ and $3\alpha+3$, respectively. Proposition~\ref{compare} and $8T_0Q_k=\Phi_{A_k^{[4]}}$ give
\[
 8T_0\psi_{2k+1}=\psi_{2k}+\epsilon+(s_k-4r_k)t .
\]
For $k\ge1$, let $J=(A_k^{[2]}:h)$. By Theorem~\ref{transitions} its cycle has half-length $2k-1$ and indices $6N_k-4$, $6N_k+1$. In $A_{k-1}^{[2]}$ the cycle has the same half-length and indices $6N_{k-1}+2$, $6N_{k-1}+4$. The two increases are $18N_{k-1}-6$ and $18N_{k-1}-3$. Now $8T_1P_k=r_k+\Phi_J$ and $\Phi_J(1)=8p_k-r_k$ by \eqref{halves}. A second application of Proposition~\ref{compare} gives
\[
 8T_1\psi_{2k}=\psi_{2k-1}+\epsilon+(r_k-4p_k)(1-t).
\]

At $k=0$ the remaining identity follows directly from $h^q=(xyz)^q$ modulo $(x^{2q},y^{2q},z^{2q})$. Multiplication by $(xyz)^q$ identifies $P/(x^q,y^q,z^q)$ with its image in the larger quotient. Hence $E_{2q}(q+j)=7q^3+E_q(j)$. Since $p_0=1$, this gives $8T_1\psi_0=\psi_0+3(1-t)$.

Together with \eqref{slopes}, these identities give the full recursion \eqref{rules}. Induction on the dyadic denominator yields $\psi_m=\phi_m$ for all $m$. Finally $A_0=\m$, so $q^{-3}E_q(j)=P_0(j/q)=j/q+\phi_0(j/q)$ for $0\le j\le q$. For $j\ge q$, the assertion follows from $h^q=0$ in $P/(x^q,y^q,z^q)$.
\end{proof}

\subsection{Rational exponents and the threshold}

Blickle, Schwede and Tucker asked whether the signature of a pair takes rational values at rational exponents \cite[Question~4.7]{BST}. They single out this cubic in \cite[Remark~4.10]{BST}. Their question is stated for the localized polynomial ring. Completion preserves the finite colengths in Proposition~4.1 of their paper, so the dyadic values agree in the two rings. Continuity gives equality for every real exponent.

\begin{corollary}\label{cor:rational-exponents}
Let $h=x^3+y^3+xyz$ and let $A$ be either $\mathbb F_2[x,y,z]_{(x,y,z)}$ or $\mathbb F_2[[x,y,z]]$. Then $s(A,h^t)$ is rational for every rational $t\ge0$. In particular, $s(A,h^{1/3})=16/49$.
\end{corollary}
\begin{proof}
Write $\mathfrak n=(x,y,z)A$. For $t=j/q\in[0,1]$, with $q=2^e$, the principal-pair formula \cite[Proposition~4.1]{BST} gives
\[
 s(A,h^t)=q^{-3}\ell(A/(\mathfrak n^{[q]}:h^j))
          =1-q^{-3}E_q(j)=1-t-\phi_0(t).
\]
The middle equality is rank--nullity for multiplication by $h^j$ on $A/\mathfrak n^{[q]}$. The signature is nonnegative and nonincreasing, and its value at $1$ is zero. It remains to treat nondyadic rationals in $(0,1)$.

Write $t=0.b_1b_2\ldots$ in base two and put $t_n=\{2^nt\}$. Starting at $m_0=0$, read one digit at a time in \eqref{rules}. At step $n$ select the row for $m_n$ and the column for $b_{n+1}$. Let $m_{n+1}$ be the next index and let $r_n$ be the remaining term in that entry, evaluated at $t_{n+1}$. Then
\[
 1-t-s(A,h^t)=\sum_{n\ge0}8^{-n-1}r_n .
\]
Indeed, for a binary expansion cut off after $N$ digits the recursion gives a finite sum, since all $\phi_m(0)$ vanish. For each fixed $n$ the truncated tails tend to $t_{n+1}$. The summands are uniformly bounded in absolute value by $(4n+4)8^{-n-1}$, since $m_n\le n$. Passing to the limit uses continuity of the signature \cite[Theorem~3.2]{BST}.

For rational $t$ the digits and tails are eventually periodic. If the index returns to zero infinitely often, two such returns occur at the same phase of the period. The subsequent indices and terms $r_n$ then repeat. Otherwise the index stays positive from some point on and changes parity at every step. Take an even multiple of the digit period. Over each successive block of this length the increments of the index repeat, so the index is affine in the block number at each phase. The table now gives the same property for $r_n$, with rational coefficients. Summing the resulting geometric series and their first moments gives a rational value in both cases.

For $t=1/3=0.\overline{01}_2$ the table gives $m_n=n$ and $r_n=2(4n+3)/3$. Thus $1-t-s(A,h^t)=\frac23\sum_{n\ge0}(4n+3)8^{-n-1}=50/147$, as claimed.
\end{proof}

Blickle, Schwede and Tucker also asked when the signature function is piecewise polynomial \cite[Question~4.8]{BST}. The first part of their question concerns a fixed characteristic. For the nodal cubic in characteristic two, the behavior at the $F$-pure threshold rules out a finite polynomial partition. The second part concerns variation of the pieces as the characteristic tends to infinity.

\begin{corollary}\label{cor:threshold}
With $A$ and $h$ as in Corollary~\ref{cor:rational-exponents}, put $S(t)=s(A,h^t)$. The $F$-pure threshold of $h$ is $1$. The function $S$ agrees with no polynomial on any interval $(1-\varepsilon,1)$ with $0<\varepsilon\le1$. Every value in $[1,29/27]$ is an accumulation value of $S(t)/(1-t)^3$ as $t\to1^-$. Using Theorem~\ref{main}, this interval extends to $[1,54/49]$.
\end{corollary}
\begin{proof}
Modulo $(x^q,y^q,z^q)$, the coefficient of $(xyz)^{q-1}$ in $h^{q-1}$ is one, whereas $h^q=0$. Thus the largest surviving exponent is $q-1$ and the threshold is $1$. At dyadic $t=j/q$ the elementary identity $E_{2q}(q+j)=7q^3+E_q(j)$ above gives $S((1+t)/2)=S(t)/8$ by the principal-pair formula. Continuity extends this to $[0,1]$.

Set $G(u)=S(1-u)$. Then $G(u/2)=G(u)/8$, so $G(u)/u^3$ is constant along each sequence $u, u/2, u/4,\ldots$. At $u=1$ this quotient equals one. At $u=3/4$, the reduced cubic formula $E_4(1)=35$ gives $S(1/4)=29/64$, so the quotient equals $29/27$. The intermediate value theorem on $[3/4,1]$ gives every value between these endpoints, and repeated halving realizes each as a limit at zero. Corollary~\ref{cor:rational-exponents} gives the further value $S(1/3)/(2/3)^3=54/49$.

If $G$ agreed with a polynomial $P$ near zero, then $P(u/2)=P(u)/8$ would force every coefficient except that of $u^3$ to vanish. The quotient would then be constant near zero, contrary to the two distinct values above.
\end{proof}

The distinction between a Frobenius sequence and an unrestricted left limit matters here. Canton's approximate left derivative uses a specified sequence approaching the threshold \cite{Canton}. In this example the normalized quotient equals one along $t=1-2^{-n}$ and has no unrestricted left limit.

\section{The invariants of the cubic suspension}\label{sec:area}

We pass from the cubic to its suspension by relating the splitting number to the Hilbert--Kunz colength at each Frobenius power. The same identity applies to the families in \cite{companion}. Let $k$ be a perfect field of characteristic $p$, let $0\ne g\in(x_1,\ldots,x_r)k[[x_1,\ldots,x_r]]$, and set
\[
 Y_g=k[[u,v,x_1,\ldots,x_r]]/(uv+g),\qquad d=r+1 .
\]

\begin{proposition}\label{prop:complement}
For $A=Y_g$ and $q=p^e$,
\begin{align}
 a_e(A)&=2q^d-\ell(A/\m^{[q]}),\label{eq:finite-complement}\\
 s(A)&=2-\eHK(A),\qquad
 S_A(w)=\frac{2}{1-p^d w}-H_A(w).\label{eq:series-complement}
\end{align}
\end{proposition}
\begin{proof}
The parameter ideal $I=(x_1,\ldots,x_r,u+v)$ has $A/I=k[[u]]/(u^2)$. Its socle is generated by $\bar u$, and $I+(u)=\m$. For a maximal Cohen--Macaulay $A$-module $M$, the Gorenstein socle formula gives its free rank as $\ell(u(M/IM))$ \cite[proof of Theorem~11(2)]{HL02}. Indeed, each free summand contributes one. If $u\bar n\ne0$ for some element of a summand $N$, the map $A/I\to N/IN$ taking $1$ to $\bar n$ is injective, since every nonzero ideal of $A/I$ contains $\bar u$. It splits because $A/I$ is Artinian Gorenstein. The retraction lifts to $N\to A$ by $\operatorname{Ext}_A^i(N,A)=0$ for $i>0$ and reduction along the regular sequence generating $I$. Its value on $n$ is a unit, so $N$ has a free summand.

Apply the formula to $M=F_*^eA$. Since $k$ is perfect, restriction of scalars by Frobenius preserves finite lengths. Thus
\[
 \ell(u(M/IM))=\ell(A/I^{[q]})-\ell(A/\m^{[q]}).
\]
The regular sequence gives $\ell(A/I^{[q]})=q^d\ell(A/I)=2q^d$. Taking limits and generating series gives \eqref{eq:series-complement}.
\end{proof}

Let $B_q=k[[x_1,\ldots,x_r]]/(x_1^q,\ldots,x_r^q)$ and put $E_q(g,j)=\dim_k B_q/(g^j)$. At $g=h$ this is the colength $E_q(j)$ of Theorem~\ref{main}.

\begin{lemma}\label{lem:boxsum}
For $q=p^e$,
\[
 \ell(Y_g/\m^{[q]})=q^r+2\sum_{j=1}^{q-1}E_q(g,j),\qquad
 a_e(Y_g)=q^r+2\sum_{j=1}^{q-1}\rk(g^j\mid B_q).
\]
If $g$ is monic of degree $b\ge1$ in one variable, then $s(Y_g)\ge1/b$ and $Y_g$ is strongly $F$-regular.
\end{lemma}
\begin{proof}
In $k[u,v]/(u^q,v^q)$, multiplication by $-uv$ has one chain of length $q$, starting at $1$, and two of length $q-i$, starting at $u^i$ and $v^i$, for $1\le i<q$. Imposing $g=-uv$ gives $B_q/(g^j)$ for a chain of length $j$. Since $g^q=0$ in $B_q$, this gives the first formula \cite[Theorem~2.2]{Mon08}, \cite[Lemma~2.4(1)]{ASS}. The second follows from \eqref{eq:finite-complement} and rank--nullity.

Suppose $g$ is monic in $t$. If $bj<q$, multiplication by $g^j$ is injective on the span of the monomials of $t$-degree less than $q-bj$. Indeed, the highest power of $t$ in a nonzero product survives reduction modulo $t^q$. Thus the rank is at least $q^{r-1}\max(q-bj,0)$. Dividing the second formula by $q^{r+1}$ and taking a Riemann sum gives $s(Y_g)\ge2\int_0^{1/b}(1-bx)\,dx=1/b$. Strong $F$-regularity follows from positivity of the signature \cite{AL03}.
\end{proof}

\subsection{The quadratic series}\label{sec:numbers}

In characteristic two, following Monsky \cite{Mtrans}, set
\begin{equation}\label{eq:theta-definition}
 \theta_g(w)=(1-2^{r+1}w)H_{Y_g}(w),\qquad
 \theta_g(2^{-r-1})=\eHK(Y_g).
\end{equation}
The second identity follows by telescoping. At $w=2^{-r-1}$, the partial sum of $\theta_g$ through degree $e$ is $2^{-e(r+1)}\ell(Y_g/\m^{[2^e]})$.

\begin{proof}[Proof of Corollaries~\ref{main:irr} and~\ref{main:series}]
For the cubic $h$ put
\[
 D=1-6w+8w^2,\qquad B=1-2w,\qquad C=2w-8w^2-24w^3 .
\]
Monsky's calculation from the colength formula of Theorem~\ref{main} gives
\begin{equation}\label{eq:theta}
 D\theta_h=B+\frac{C}{\sqrt{1-4w^2}},
\end{equation}
where the square root has constant term $1$ \cite[Lemma~2.1]{Mtrans}. At $w=1/16$ the three polynomials have values $21/32$, $7/8$ and $45/512$, and the square root is $3\sqrt7/8$. Therefore $\eHK(R)=\theta_h(1/16)=4/3+5/(14\sqrt7)$. The signature follows from \eqref{eq:series-complement}.

For the series, put $L=(1-16w)D$ and $P=2D-B=1-10w+16w^2$. Since $LS_R=P-C/\sqrt{1-4w^2}$, we obtain
\begin{equation}\label{eq:splitting-polynomial}
 (1-4w^2)\bigl(L(w)S_R(w)-P(w)\bigr)^2=C(w)^2 .
\end{equation}
Each of $H_R$ and $S_R$ recovers the square root, since $C\ne0$. The radicand has simple zeros at $\pm1/2$. This proves that each series generates the stated quadratic extension.
\end{proof}

Corollary~\ref{main:series} settles the conditional part of Caminata and Zerman's example \cite[Example~4.6]{CZ}. They show that the colength function $\phi_{h,2}(j/q)=q^{-3}E_q(j)$ is a weak $2$-fractal. If it were a $2$-fractal, their argument would give a rational Hilbert--Kunz series for $uv+h$. Corollary~\ref{main:series} rules this out.

Monsky's algebraicity theorem also applies after adding a strongly rational series in disjoint variables. For a finite field $k$ of characteristic two, such a series $g$ in $r$ variables is characterized by the finite-dimensional span over $\mathbb Q$ of the iterates of $T_0,T_1$ on $j/q\mapsto q^{-r}E_q(g,j)$ \cite{MT,M09}. Since $uv+g$ is again strongly rational, Theorem~\ref{main} and \cite[Theorem~3.4]{M09} give algebraic Hilbert--Kunz series and multiplicity for $Y_{h+g}$. Equation~\eqref{eq:series-complement} gives algebraicity of its splitting series and signature. This includes $g=\sum_i g_i(x_i,y_i)$ in disjoint pairs \cite[Theorem~3.3]{M09}.

\section{An elliptic period and six hypersurfaces}\label{sec:period}

Monsky's product formula expresses the series of a sum in disjoint variables as a coefficientwise product. With $\odot$ denoting coefficientwise multiplication, his formula reads \cite[Theorem~1.6]{Mtrans}
\begin{equation}\label{eq:hadamard}
 \theta_{f_1+\cdots+f_k}=\theta_{f_1}\odot\cdots\odot\theta_{f_k}.
\end{equation}
To isolate the even and odd coefficients we also use the following identities, where $f$ involves $r$ variables \cite[Theorem~1.8 and Lemma~1.9]{Mtrans}.
\begin{align}
 \theta_{f^2}&=1+2^r w\theta_f,\label{eq:square-theta}\\
 \theta_{t^5}-\theta_{t^3}&=2\sum_{n\ge1}w^{2n},\qquad
 2\theta_{t^3}-\theta_{t^5}=1+2\sum_{n\ge0}w^{2n+1}.
 \label{eq:parities}
\end{align}

\begin{proof}[Proof of Theorem~\ref{main:trans}]
Write $\mu_{i,b}=\eHK(A_{i,b})$. The six hypersurfaces have dimension eight and are Gorenstein. Each sum $h(X)+h(X')^{2^i}+t^b$ is monic of degree three in the first coordinate of $X$. Lemma~\ref{lem:boxsum} gives $s(A_{i,b})\ge1/3$ and strong $F$-regularity for all six rings.

Consider the number
\[
 T=\sum_{n\ge0}\frac{\binom{2n}{n}^2}{65536^n}.
\]
Monsky relates this number to Hilbert--Kunz multiplicities \cite[Theorems~1.10 and~2.2]{Mtrans}. For the six rings in Theorem~\ref{main:trans}, the relation is
\begin{equation}\label{eq:six-multiplicities}
 \begin{split}
 2T+94={}&-737\mu_{0,3}+881\mu_{0,5}
        +1440\mu_{1,3}-1488\mu_{1,5}\\
        &{}-256\mu_{2,3}+256\mu_{2,5}.
 \end{split}
\end{equation}

To compute the coefficients, let $b_n=[w^n]\theta_h$. Comparing the coefficients of $w^{2n+2}$ in \eqref{eq:theta} gives
\[
 b_{2n+2}-6b_{2n+1}+8b_{2n}=-8\binom{2n}{n},\qquad b_0=1 .
\]
Put $\xi=1/256$, $z=\xi^2$, and write $\langle d_n\rangle=\sum_{n\ge0}d_n z^n$. By \eqref{eq:theta-definition} and \eqref{eq:hadamard}, $\mu_{i,b}$ is the value of $\theta_h\odot\theta_{h^{2^i}}\odot\theta_{t^b}$ at $\xi$. Write $U_i=\mu_{i,5}-\mu_{i,3}$ and $V_i=2\mu_{i,3}-\mu_{i,5}$. Equations~\eqref{eq:square-theta} and \eqref{eq:parities} give the six terms in the square of the recurrence.
\[
\renewcommand{\arraystretch}{1.25}
\begin{array}{l@{\ }c@{\ }l@{\qquad}l@{\ }c@{\ }l}
 \langle b_{2n}^2\rangle&=&1+U_0/2&
 \langle b_{2n+2}^2\rangle&=&U_0/(2z),\\
 \langle b_{2n+1}^2\rangle&=&(V_0-1)/(2\xi)&
 \langle b_{2n+1}b_{2n+2}\rangle&=&U_1/(16z),\\
 \langle b_{2n}b_{2n+1}\rangle&=&(V_1-1)/(16\xi)&
 \langle b_{2n}b_{2n+2}\rangle&=&U_2/(128z).
\end{array}
\]
For instance, $[w^{2n+2}]\theta_{h^2}=8b_{2n+1}$ gives $U_1=16z\langle b_{2n+1}b_{2n+2}\rangle$. Substitution into $64T=\langle(b_{2n+2}-6b_{2n+1}+8b_{2n})^2\rangle$ proves \eqref{eq:six-multiplicities}.

Euler's integral identifies $2\pi T$ with a nonzero real period of the smooth projective model of $y^2=(1-x^2)(1-16\xi^2x^2)$. By Schneider's theorem this period divided by $\pi$ is transcendental \cite{Schneider}, \cite[Lemma~2.4 and Theorem~2.5]{Mtrans}. Equation~\eqref{eq:six-multiplicities} then forces at least one $\mu_{i,b}$ to be transcendental. For the same pair, $s(A_{i,b})=2-\mu_{i,b}$ by Proposition~\ref{prop:complement}.
\end{proof}

\section{Sums of nodal cubics}\label{sec:repeated}

Monsky and Teixeira proved rationality of the Hilbert--Kunz series for sums of bivariate series in disjoint pairs of variables over a finite field \cite{MTII}. Here each summand is a cubic in three variables. For $g_m=\sum_j h(X_j)$, the product formula \eqref{eq:hadamard} raises each coefficient $b_n=[w^n]\theta_h$ to the $m$th power. Thus
\begin{equation}\label{eq:repeated}
 (1-2^{3m+1}w)H_{R_m}(w)=\sum_{n\ge0}b_n^m w^n .
\end{equation}
The second term in the asymptotic expansion of $b_n$ determines transcendence for $m\ge2$.

\begin{lemma}\label{lem:secondterm}
For $\alpha=1-\sqrt3/2$ and $c_n=b_n-\alpha4^n$, one has
\[
 c_n\sim 4\sqrt{2/\pi}\,2^n\sqrt n .
\]
\end{lemma}
\begin{proof}
Set $J(w)=(1-4w)\theta_h(w)=\sum d_n w^n$. Equation~\eqref{eq:theta} gives
\[
 J(w)=1+\frac{P(w)}{(1-4w^2)^{3/2}},\qquad
 P(w)=2w(1-6w)(1+2w)^2 .
\]
In particular $J(1/4)=\alpha$. The coefficient of $w^{2k}$ in $(1-4w^2)^{-3/2}$ is $(2k+1)\binom{2k}{k}$, and the odd coefficients vanish. By Stirling's formula the ratio of this even coefficient to $2^{2k}\sqrt{2k}$ tends to $\sqrt{2/\pi}$. In $d_n/(2^n\sqrt n)$ the remaining factor is therefore the even or the odd part of $P$ evaluated at $1/2$, according to the parity of $n$. Both have value $-4$, since $P(1/2)=-8$ and $P(-1/2)=0$.

Expansion of $J/(1-4w)$ gives $b_n=4^n\sum_{i=0}^n d_i4^{-i}$. Subtract $4^n J(1/4)$ to obtain $c_n=-\sum_{j\ge1}4^{-j}d_{n+j}$. After division by $2^n\sqrt n$, the $j$th term is bounded by a constant times $2^{-j}\sqrt{1+j}$ for $n\ge1$. Dominated convergence and $\sum_{j\ge1}2^{-j}=1$ give the stated coefficient.
\end{proof}

\begin{proof}[Proof of Theorem~\ref{main:repeated}]
Fix $m\ge2$ and put $\lambda=4^{m-1}$. The series $U(w)=\sum c_n w^n=\theta_h(w)-\alpha/(1-4w)$ is algebraic over $\overline{\mathbb Q}(w)$. Remove the terms constant and linear in $c_n$ from $b_n^m=(\alpha4^n+c_n)^m$, and write the remainder as
\[
 d_{n,m}=b_n^m-\alpha^m4^{mn}-m\alpha^{m-1}\lambda^n c_n .
\]
Since $c_n/4^n\to0$, the quadratic term determines the limit. Lemma~\ref{lem:secondterm} gives
\[
 \frac{d_{n,m}}{n\lambda^n}\longrightarrow
 K_m=\binom m2\alpha^{m-2}\frac{32}{\pi}.
\]
For $G_m(w)=\sum_{n\ge0}d_{n,m}w^n$, it follows that
\[
 \lim_{w\uparrow 1/\lambda}(1-\lambda w)^2G_m(w)=K_m .
\]
Indeed, for every $\varepsilon>0$ all sufficiently large coefficients differ from $K_m n\lambda^n$ by at most $\varepsilon n\lambda^n$. Summing this estimate uses $\sum_{n\ge0}n x^n=x/(1-x)^2$.

An algebraic function over $\overline{\mathbb Q}(w)$ with a finite limit at an algebraic point has algebraic limit along any real interval on which it is defined. To see this, clear denominators in a polynomial relation and divide its coefficients by their common power of the local parameter before specializing. Here $K_m$ is a nonzero algebraic multiple of $1/\pi$, hence transcendental. Thus $G_m$ is transcendental. The two series removed from the right side of \eqref{eq:repeated} are $\alpha^m/(1-4^m w)$ and $m\alpha^{m-1}U(\lambda w)$, both algebraic. This proves transcendence of $H_{R_m}$, and \eqref{eq:series-complement} gives it for $S_{R_m}$.

Every algebraic series in characteristic zero satisfies a linear differential equation with polynomial coefficients. This property is preserved by Hadamard products \cite[Theorem~2.10]{Stanley}, as well as by addition and multiplication by rational functions. Equations~\eqref{eq:repeated} and \eqref{eq:series-complement} give the differential equations asserted for both series. Finally $g_m$ is monic of degree three in the first coordinate of $X_1$. Lemma~\ref{lem:boxsum} gives the lower bound for the signature and strong $F$-regularity.
\end{proof}

\begin{corollary}
For every $m\ge1$, the ring $R_m$ has infinite $F$-representation type.
\end{corollary}
\begin{proof}
Suppose that only finitely many isomorphism classes of indecomposable modules occur among the $F_*^eR_m$, and choose representatives $M_1,\ldots,M_r$. Since $R_m$ is complete, Krull--Schmidt gives well-defined multiplicities. The pushforward of each $M_i$ is a direct summand of some $F_*^{e+1}R_m$. Thus a fixed integer matrix $B$ takes the multiplicity vector at level $e$ to that at level $e+1$. Let $v$ be the initial vector and let $c_i$ be the minimal number of generators of $M_i$. Since the residue field is $\mathbb F_2$, reduction of $F_*^eR_m$ modulo the maximal ideal gives a vector space of dimension $\ell(R_m/\mathfrak m^{[2^e]})$. Hence
\[
 H_{R_m}(w)=c^{\mathsf T}(I-wB)^{-1}v.
\]
This contradicts Corollary~\ref{main:series} for $m=1$ and Theorem~\ref{main:repeated} for $m\ge2$.
\end{proof}

The two series in Theorem~\ref{main:repeated} are diagonals of rational functions in $2m$ variables. By the algebraic diagonal theorem \cite{DL}, write $\theta_h=\operatorname{Diag} F(x,y)$ with $F$ rational and regular at the origin. Put $U=\prod_{i=1}^m x_i y_i$. Extracting the coefficient with all exponents equal in \eqref{eq:repeated} gives
\[
 H_{R_m}(w)=\operatorname{Diag}\left(
 \frac{\prod_{i=1}^m F(x_i,y_i)}{1-2^{3m+1}U}\right).
\]
The same assertion for $S_{R_m}$ follows from \eqref{eq:series-complement}. Thus this family satisfies the conclusion of Christol's conjecture on globally bounded $D$-finite series \cite{Christol}.

In the terminology of Allouche and Mend\`es France \cite{AMF}, the Hadamard grade over $\mathbb Q$ is the least number of algebraic series whose coefficientwise product is the given series. By \eqref{eq:repeated}, $(1-2^{3m+1}w)H_{R_m}(w)$ has grade at most $m$, since it is the $m$th Hadamard power of $\theta_h$. For $m=2$ its grade is exactly two by Theorem~\ref{main:repeated}.

\section{The multiplicity for two nodal cubics}\label{sec:r2}

For $R_2$, equation~\eqref{eq:repeated} reduces the multiplicity to a value of the Hadamard square of $\theta_h$. We express this value as a period and determine its arithmetic nature.
\begin{lemma}\label{lem:period}
Let
\[
 a=\frac1{128},\qquad
 Q(z)=\frac{2z(1-6z)(1+2z)}{(1-2z)(1-4z)},\qquad
 T(z)=\frac{Q(z)}{\sqrt{1-4z^2}}.
\]
Then
\begin{equation}\label{eq:split}
 e_{\mathrm{HK}}(R_2)=\frac87+\frac{221}{105\sqrt{255}}+I,
 \qquad I=(T\odot T)(a).
\end{equation}
Moreover, on the smooth genus-one curve
\begin{equation}\label{eq:curve}
 E:\quad y^2=(1-4z^2)(z^2-4a^2)
       =(1-4z^2)\left(z^2-\frac1{4096}\right),
\end{equation}
put $\Omega=Q(z)Q(a/z)\,dz/y$. There is a closed cycle $\gamma$ on
$E(\mathbb C)$ such that
\begin{equation}\label{eq:period}
 I=\frac{1}{2\pi i}\int_\gamma\Omega,
 \qquad \int_\gamma\frac{dz}{y}
  =2\pi i\sum_{n\geq0}\binom{2n}{n}^{\!2}a^{2n}\neq0.
\end{equation}
In particular, the image of $\gamma$ in $H_1(E,\mathbb Z)$ is nonzero.
\end{lemma}
\begin{proof}
Equation~\eqref{eq:theta} simplifies to
$\theta_h(z)=(1-4z)^{-1}+T(z)$. The Hadamard product of two copies at $a$ is
\[
 (\theta_h\odot\theta_h)(a)=\frac1{1-16a}+2T(4a)+(T\odot T)(a).
\]
Since $4a=1/32$ and $Q(1/32)=221/3360$, we get
$2T(4a)=221/(105\sqrt{255})$. The equality
$e_{\mathrm{HK}}(R_2)=(\theta_h\odot\theta_h)(a)$ follows by taking the radial limit in \eqref{eq:repeated}. Indeed, $b_n=O(4^n)$ by Lemma~\ref{lem:secondterm}, so the Hadamard square converges at $a=1/128$. On the left, put $w=x/128$. The coefficients of $H_{R_2}(x/128)$ tend to $\eHK(R_2)$ by definition, so Abel's limit theorem gives the same value as $x\uparrow1$ after multiplication by $1-x$.

For $2a<r<1/4$ and $a/r<1/4$ (in particular $r=\sqrt a$), Cauchy's coefficient formula gives
\[
 I=\frac1{2\pi i}\oint_{|z|=r}T(z)T(a/z)\,\frac{dz}{z}.
\]
Choose the branch on this circle given by
$y=z\sqrt{1-4z^2}\sqrt{1-4a^2/z^2}$. The integral is
$\frac1{2\pi i}\int_\gamma\Omega$. The same expansion with $Q\equiv1$ yields
\[
 \frac1{2\pi i}\int_\gamma\frac{dz}{y}
 =\sum_{n\geq0}\binom{2n}{n}^2a^{2n}>0.
\]
Since $dz/y$ is a holomorphic differential on $E$, the homology class of $\gamma$ cannot vanish.
\end{proof}

\begin{lemma}\label{lem:divisor}
Let
\[
 P_\pm=\left(\frac14,\,\pm\frac{3\sqrt{85}}{128}\right),\qquad
 Q_\pm=\left(\frac1{32},\,\pm\frac{3\sqrt{85}}{1024}\right)
\]
be points on $E$. The nonzero residues of $\Omega$ are
\[
 \res_{P_+}\Omega=\rho,\quad
 \res_{P_-}\Omega=-\rho,\quad
 \res_{Q_+}\Omega=-\rho,\quad
 \res_{Q_-}\Omega=\rho,
 \qquad \rho=\frac{221}{420\sqrt{85}}\neq0.
\]
All other poles of $\Omega$ are of the second kind. Let
\[
 \mathcal D=(P_+)+(Q_-)-(P_-)-(Q_+).
\]
Take $O=(2a,0)$ as the origin of $E$ and $Z\in E(\overline{\mathbb Q})$ with
$[\mathcal D]=[(Z)-(O)]$ in $\Pic^0(E)$. Then $Z=4P_+$ in the group law and $Z$ has infinite order.
\end{lemma}
\begin{proof}
The rational function $Q(z)$ has residue $3/16$ at $z=1/4$; the function
$Q(a/z)$ has residue $-3a$ at $z=4a=1/32$.
Using $Q(1/32)=221/3360$ and the displayed $y$-coordinates gives the four residues. The poles at $z=1/2$ and $z=2a$ lie over branch points. A local parameter $t$ with $z-z_0=t^2$ shows that their principal parts have no residues. The putative poles at $z=0,\infty$ are removable.

The hyperelliptic involution $\sigma(z,y)=(z,-y)$ is negation on the elliptic curve with origin $O$. The automorphism
\[
 \iota(z,y)=\left(\frac{a}{z},\frac{ay}{z^2}\right)
\]
satisfies $\iota^*(dz/y)=-dz/y$ and sends $O$ to the two-torsion point
$U=(1/2,0)$. Hence $\iota(R)=U-R$ for $R\in E$. Since
$\iota(P_+)=Q_+$, we have
$P_-=-P_+$, $Q_+=U-P_+$, and $Q_-=P_+-U$. Thus the group sum of $\mathcal D$ is
$4P_+-2U=4P_+$.

To show that $P_+$ is non-torsion, use the degree-two map
\[
 \pi:E\longrightarrow E',\qquad (z,y)\longmapsto
 (X,Y)=(-4z^2,4zy),
\]
where
\begin{equation}\label{eq:weierstrass}
 E':\quad Y^2=X(X+1)\left(X+\frac1{1024}\right).
\end{equation}
The point $P_+$ maps to
$P'=(-1/4,3\sqrt{85}/128)\in E'(\mathbb Q(\sqrt{85}))$.
Both primes $7$ and $19$ split in $\mathbb Q(\sqrt{85})$, and $E'$ has good reduction at these primes. Choosing $\sqrt{85}\equiv1\pmod7$ and $\sqrt{85}\equiv3\pmod{19}$, the reductions are
\[
 \begin{array}{c|c|c|c}
 \text{prime}&\overline{P'}&2\overline{P'}&4\overline{P'}\\\hline
 7&(5,5)&(0,0)&O_{E'}\\
 19&(14,2)&(16,6)&(0,0)
 \end{array}
\]
so their respective orders are $4$ and $8$. Since good reduction in odd residue characteristic is injective on $2$-power torsion, $P'$ cannot be torsion. The morphism $\pi$, after translating by $-\pi(O)$, is an isogeny. Since $\pi(O)=(-1/1024,0)$ is torsion, $P_+$ is also non-torsion.
\end{proof}

\begin{lemma}\label{lem:nocm}
The elliptic curve $E$ has no complex multiplication.
\end{lemma}
\begin{proof}
The isogenous curve $E'$ in \eqref{eq:weierstrass} has Legendre parameter $\lambda=1/1024$ and
\[
 j(E')=256\frac{(1-\lambda+\lambda^2)^3}{\lambda^2(1-\lambda)^2}.
\]
This rational number is not an integer, since $v_3(1-\lambda)>0$ but
$v_3(1-\lambda+\lambda^2)=0$. A CM $j$-invariant is an algebraic integer. Thus $E'$ has no CM, and the same holds for its isogenous curve $E$.
\end{proof}

\begin{proof}[Proof of Theorem \ref{main:r2}]
By Lemma \ref{lem:divisor}, the divisor $\mathcal D$ represents the non-torsion point $Z=4P_+$ of $\Pic^0(E)$. Hence there is a rational function $f$ on $E$, defined over $\overline{\mathbb Q}$, with
\[
 \divi(f)=\mathcal D-(Z)+(O).
\]
Choose a third-kind differential $\xi_Z$ defined over $\overline{\mathbb Q}$ with residues $+1$ at $Z$ and $-1$ at $O$ and no other poles. Any two such choices differ by a holomorphic differential, so the choice changes its periods only by first-kind periods. By Lemma \ref{lem:divisor}, the differential
\[
 \Omega-\rho\,\xi_Z-\rho\,d\log f
\]
has no residues, hence is of the second kind. On a smooth projective curve in characteristic zero, differentials of the second kind modulo exact differentials represent its first de Rham cohomology \cite{Mukai,HuberWuestholz2022}. Hence
\begin{equation}\label{eq:decompose}
 \Omega=\rho\,\xi_Z+\rho\,d\log f
      +A\,\omega_E+B\,\eta_E+dF,
 \qquad A,B\in\overline{\mathbb Q},
\end{equation}
where $\omega_E$ and $\eta_E$ represent a basis of
$H^1_{\mathrm{dR}}(E/\overline{\mathbb Q})$ and $F$ is a rational function on $E$.

Remove the poles of all the differentials in \eqref{eq:decompose} from $E$. A small deformation of $\gamma$ avoids the removed points and preserves its $\Omega$-period. Choose lifts $\gamma_1,\gamma_2$ of an integral homology basis to the punctured curve. By Lemma~\ref{lem:period}, the image of $\gamma$ on the compact curve has coordinates $(m,n)\ne(0,0)$. On the punctured curve the difference $\gamma-m\gamma_1-n\gamma_2$ is a sum of small loops about the punctures. Each contributes $2\pi i$ times an algebraic residue. Also, $d\log f$ has periods in $2\pi i\mathbb Z$. Thus, modulo the $\overline{\mathbb Q}$-span of $2\pi i$ and the periods of $\omega_E,\eta_E$, we have
\[
 \int_\gamma\Omega\equiv\rho\left(m\int_{\gamma_1}\xi_Z+n\int_{\gamma_2}\xi_Z\right).
\]

Write $\omega_i,\eta_i,\xi_i$ for the integrals of $\omega_E,\eta_E,\xi_Z$ along $\gamma_i$. The elliptic period theorem \cite[Theorem~2 and Section~3.4]{Wuestholz2021} gives dimension eight for the period space over $\overline{\mathbb Q}$ when the elliptic curve has no complex multiplication and its point in $\Pic^0$ has infinite order. Both hypotheses hold by Lemmas~\ref{lem:nocm} and~\ref{lem:divisor}. The eight generators
\[
 1,\quad 2\pi i,\quad \omega_1,\omega_2,\quad
 \eta_1,\eta_2,\quad \xi_1,\xi_2
\]
are therefore independent. Changing a lift of a cycle adds a multiple of $2\pi i$ to its third-kind period. Changing the third-kind differential with the same simple poles and residues adds first-kind periods. A change of algebraic de Rham basis acts by an invertible matrix over $\overline{\mathbb Q}$ on each pair $(\omega_i,\eta_i)$; exact differentials have zero period. These changes preserve independence. Since $\rho\ne0$ and $(m,n)\ne(0,0)$, the preceding congruence gives
$\int_\gamma\Omega\notin2\pi i\,\overline{\mathbb Q}$, and
$I=(2\pi i)^{-1}\int_\gamma\Omega$ is transcendental. Equation \eqref{eq:split} shows that $e_{\mathrm{HK}}(R_2)$ is transcendental. Finally
$s(R_2)=2-e_{\mathrm{HK}}(R_2)$ is transcendental as well.
\end{proof}

\begin{corollary}
The $F$-rational signature of Hochster and Yao \cite{HY}, the relative $F$-rational signature of Smirnov and Tucker, and the dual $F$-signature of Sannai \cite{Sannai} are all transcendental for $R_2$. Every positive-degree Frobenius Betti number of its residue field is transcendental as well.
\end{corollary}
\begin{proof}
The hypersurface $R_2$ is Gorenstein. Its $F$-rational and relative $F$-rational signatures therefore equal $s(R_2)$ by the socle formula \cite[Introduction and Proposition~3.15]{ST}. Its canonical module is free of rank one, so every surjection $F_*^e\omega_{R_2}\to\omega_{R_2}^{\oplus n}$ splits. The largest such $n$ is the free rank of $F_*^eR_2$. Taking limits gives the same equality for the dual signature. For the residue field $k=\mathbb F_2$, the hypersurface formula of De Stefani, Huneke and N\'u\~nez-Betancourt \cite[Example~3.2]{DSHN} gives, for every $i\ge1$ and $m\ge1$,
\[
 \beta_i^F(k,R_m)=\eHK(R_m)-s(R_m)=2\eHK(R_m)-2.
\]
For $m=1$ this is $2/3+5/(7\sqrt7)$ by Corollary~\ref{main:series}. For $m=2$ the assertions follow from Theorem~\ref{main:r2}.
\end{proof}

\begin{corollary}\label{cor:monotone}
For $m\ge1$ put $f_m=uv+\sum_{j=1}^m h(X_j)$ and let $C_{m,r}$ be the complete hypersurface over $\mathbb F_2$ defined by $f_m+z^{2^r}$, with $z$ a new variable. Write $E_{m,r}=\eHK(C_{m,r})$, $b_r=[w^r]\theta_h(w)$ and $\alpha=1-\sqrt3/2$. Then $E_{m,r}$ is dyadic rational and
\[
 E_{m,r}-E_{m,r-1}=\frac{b_r^m}{2^{(3m+1)r}}\quad(r\ge1),\qquad
 \eHK(R_m)-E_{m,r}\sim
 \frac{\alpha^m}{2^{m+1}-1}\,2^{-(m+1)r}\quad(r\to\infty).
\]
For each fixed $m$ the multiplicities eventually increase strictly and the signatures eventually decrease strictly. For even $m$ the multiplicities are nondecreasing from the start.
\end{corollary}
\begin{proof}
Let $A=\mathbb F_2[[x_1,\ldots,x_n]]$, let $0\ne f\in\mathfrak m_A$, and set $B=A/(f)$. For $Q=2^r$ put $C_Q=A[[z]]/(z^Q-f)$. Aberbach and Enescu used radical extensions to bound Hilbert--Kunz multiplicities \cite{AE}, and Aberbach and Cox-Steib used them to interpolate between multiplicities \cite{ACS}. Here the Frobenius colengths can be read off exactly. For $q=QP$, with $P$ a power of two, we have $z^q=f^P$ in $C_Q$. Since $C_Q$ is free of rank $Q$ over $A$, flatness of Frobenius on $A$ gives
\[
 \ell(C_Q/\mathfrak m_{C_Q}^{[q]})
 =Q\,\ell(A/(x_1^q,\ldots,x_n^q,f^P))
 =QP^n\,\ell(B/\mathfrak m_B^{[Q]}).
\]
Divide by $q^n$ and let $P$ tend to infinity. This gives
\[
 \eHK(C_Q)=Q^{1-n}\ell(B/\mathfrak m_B^{[Q]}).
\]
The right side is dyadic rational and tends to $\eHK(B)$ as $r$ tends to infinity. For $f=f_m$, convergence also follows from the continuity theorem of Polstra and Smirnov \cite{PS}. To obtain the rate, take $n=3m+2$ and compare coefficients in \eqref{eq:repeated}. This gives the stated difference. By Lemma~\ref{lem:secondterm}, $b_r\sim\alpha4^r$, so the difference is eventually positive and asymptotic to $\alpha^m2^{-(m+1)r}$. Summing the geometric tail gives the claimed equivalent. For even $m$, every $b_r^m$ is nonnegative. Finally $s(C_{m,r})=2-E_{m,r}$ by Proposition~\ref{prop:complement}.
\end{proof}

For $m=2$ the multiplicities increase to $\eHK(R_2)$ with error asymptotic to $(\alpha^2/7)8^{-r}$. The signatures decrease to $s(R_2)\ge1/3$, so every $C_{2,r}$ is strongly $F$-regular \cite{AL03}. Thus in dimension eight over $\mathbb F_2$ there is a strictly increasing chain of dyadic Hilbert--Kunz multiplicities and a strictly decreasing chain of dyadic $F$-signatures, both with transcendental limits. Adams, Sandstrom and Simpson study the ascending chain condition for signatures \cite{ASS}. For $m=1$ the corresponding multiplicities approach $4/3+5/(14\sqrt7)$ in dimension five.

\end{document}